\documentclass[11pt]{amsart}

\usepackage{amsmath,amssymb,amscd}

\usepackage{amsrefs}
\usepackage{enumerate}
\usepackage{xcolor}

\usepackage{amsmath,amssymb}

\usepackage{amsxtra, amsmath}
\usepackage{amssymb, amscd}
\usepackage{graphicx}
\usepackage{mathrsfs}

\usepackage{mathtools}
\usepackage{float}

\newtheorem{thm}{Theorem}[section]
\newtheorem{prop}[thm]{Proposition}

\newtheorem{cor}[thm]{Corollary}

\theoremstyle{definition}
\newtheorem{definition}[thm]{Definition}

\theoremstyle{remark}
\newtheorem{remark}[thm]{Remark}

\usepackage{dsfont}
\allowdisplaybreaks
\title{Darboux transformations and the algebra $\mathcal{D}(W)$.}
\author{Ignacio Bono Parisi}
\author{Ines Pacharoni}

\subjclass[2020]{33C45, 42C05, 34L05, 34L10}
\thanks{This paper was partially supported by SeCyT-UNC, CONICET, PIP 1220150100356.}
\keywords{ Matrix-valued orthogonal polynomials, matrix Bochner problem, Darboux transformations, discrete-continuous bispectrality, matrix-valued bispectral functions}

\address{CIEM-FaMAF\\ Universidad Nacional de C\'or\-do\-ba\\
CP 5000, C\'or\-do\-ba,  Argentina}
\email{ignacio.bono@unc.edu.ar, ines.pacharoni@unc.edu.ar}

\begin{document}
\begin{abstract}

The problem of finding weight matrices $W(x)$ of size $N \times N$ such that the associated sequence of matrix-valued orthogonal polynomials are eigenfunctions of a second-order matrix differential operator is known as the Matrix Bochner Problem, and it is closely related to Darboux transformations of some differential operators. 

This paper aims to study Darboux transformations between weight matrices and to establish a direct connection with the structure of the algebra $\mathcal D(W)$ of all differential operators that have a sequence of matrix-valued orthogonal polynomials with respect to $W$ as eigenfunctions. 
\end{abstract} 
\maketitle

\section{Introduction}

The theory of matrix-valued orthogonal polynomials starts with the work of Krein, see \cites{K49, K71}. If one is considering possible applications of these matrix polynomials, it is natural to concentrate on those cases where some extra property holds. A. Durán in \cite{D97}, posed the problem of characterizing those positive definite matrix-valued weights whose matrix-valued orthogonal polynomials are common eigenfunctions of some symmetric second-order differential operator with matrix coefficients. 
This was an important generalization to the matrix-valued case of the problem originally considered by S. Bochner in 1929, see \cite{B29}.
S. Bochner proved that, up to an affine change of coordinates, the only weights in the real line  satisfying these properties are the classical weights $w(x)=e^{-x^2}$, $w_{\alpha}(x)=x^\alpha e^{-x} $ 
and 
$w_{\alpha,\beta}(x)=(1-x)^\alpha(1 + x)^\beta$ of Hermite, Laguerre, and Jacobi respectively.

The problem of finding weight matrices $W(x)$ of size $N \times N$  such that the associated sequence of orthogonal matrix polynomials are eigenfunctions of a second-order differential operator is known as the Matrix Bochner Problem (MBP). 
  In the matrix-valued case,  the first nontrivial solutions of this problem were given in  \cite{GPT01} and \cite{GPT02}, as a consequence of the study of matrix-valued spherical functions and also in  \cite{DG04} and \cite{G03}. 
In the past twenty years, a large amount of examples have been found, not necessarily associated with Lie theory.
See \cite{GPT03a}, \cite{GPT03b}, \cite{DG04}, \cite{GPT05}, \cite{CG06}, \cite{PT07}, \cite{DG07}, \cite{P08}, \cite{PR08}, \cite{PZ16}, \cite{KPR12}. 
However, the explicit construction of such examples is, even nowadays, not an easy task. 

After the appearance of the first examples, some works focused on studying the algebra $\mathcal{D}(W)$ of all differential operators that have a sequence of matrix-valued orthogonal polynomials with respect to $W$ as eigenfunctions, see \cite{T11}, \cite{Z17}.

In \cite{CY18}, W. R. Casper and M. Yakimov made a significant breakthrough toward solving the MBP. They proved that under certain assumptions on the algebra $\mathcal{D}(W)$, the solutions to the MBP can be obtained through a Darboux transformation of some direct sum of classical scalar weights.

However, this result does not fully resolve the Matrix Bochner Problem since many examples fall outside this characterization (see \cite{BP23} and \cite{BP24-1}). Nevertheless, it highlights the importance of considering Darboux transformations of diagonal weights to obtain new solutions to the MBP. 

 We say that a weight $\widetilde W$ is a Darboux transformation of a weight $W$ if there exists a differential operator $D\in \mathcal{D}(W)$ that can be factorized as $D = \mathcal{V}\mathcal{N}$, for some degree-preserving operators $\mathcal{V}, \, \mathcal{N} $ and such that the polynomials $P_n(x) \cdot \mathcal{V}  $ are orthogonal with respect to the weight $\widetilde W$, for $P_n(x)$ a sequence of orthogonal polynomials for $W$. 

To obtain new nontrivial solutions to the MBP through Darboux transformations of $W$, a direct sum of classical scalar weights, we need to factorize operators from the algebra $\mathcal{D}(W)$. For this reason, it is crucial to fully understand the entire algebra $\mathcal{D}(W)$. Many nontrivial solutions arise from the factorization of operators that are non-diagonal and of order higher than two.

In the classical scalar cases of  Hermite, Laguerre, or Jacobi weights, the structure of the algebra $\mathcal D(w)$ is well understood:
it is a polynomial algebra in the classical differential operator of Hermite, Laguerre, or Jacobi, respectively.
However, in the matrix cases, the structure of this algebra is certainly a non-trivial problem, even for diagonal weights. 

\smallskip
This paper aims to study Darboux transformations of weight matrices and their connection with the structure of the algebra $\mathcal{D}(W)$.

The paper is organized as follows.  
Section \ref{sec-MOP} contains background material about matrix-valued orthogonal polynomials and the algebra $\mathcal D(W)$.    Darboux transformations of weight matrices are studied in Section \ref{Darboux-sect}. We prove that the Darboux transformation defines an equivalence relation in the set of weight matrices defined on the same interval. We also establish some important relations between the algebras $\mathcal D(W)$ and $\mathcal D(\widetilde W)$ when $W$ is a Darboux transformation of $\widetilde W$.

In Section \ref{classical-section} and Section \ref{sub}, we focus on classical scalar Hermite, Laguerre, and Jacobi weights. The main results are the following:
\begin{itemize}
    \item A (shifted) Hermite weight 
$w_b    \text{ is not a Darboux transformation of }    w_{a} \text{  if }   a\neq b.$

\item A Laguerre weight 
 $w_\alpha  \text{ is a Darboux transformation of } w_{\widetilde \alpha}$  if and only if $\alpha-\widetilde \alpha \in \mathbb Z.$
 
\item A Jacobi weight $w_{\alpha,\beta}$   is a Darboux transformation of $w_{\widetilde{\alpha},\widetilde{\beta}}$  if and only if $\alpha +\beta= \widetilde{\alpha}+ \widetilde{\beta}$,   and  $\alpha-\widetilde{\alpha}\in \mathbb{Z}.$ 
\end{itemize}

In Section \ref{algstruct-sect}, we describe the form of any operator $D$ in the algebra $\mathcal{D}(W)$, for $W(x) = w_1(x) \oplus \cdots \oplus w_N(x)$ a direct sum of classical scalar weights of the same type. According to \cite{TZ18}, we have
\[ \mathcal{D}(W) = \sum_{i,j=1}^N \mathcal{D}(w_i, w_j) \, E_{i,j}, \]
where $\mathcal{D}(w_{i}, w_{j}) $ is the set of differential operators $\tau$ such that $\, p_{n}(x) \cdot \tau = \lambda_{n}(\tau) q_{n}(x)$,  for  all $n\in \mathbb N_0$, with $p_n$ and $q_n$  the sequence of orthogonal polynomials for the weights $w_i$ and $w_j$, respectively. As usual, $E_{i,j}$ denotes the $N \times N$ matrix whose $(i,j)$-entry is 1 and 0 elsewhere.
  In particular, for  $w_i=w_j$, we have  $\mathcal{D}(w_i,w_j) = \mathcal{D}(w_j)$.

We prove that $\mathcal D(w_i,w_j)$ is a  cyclic right-module over the algebra $\mathcal{D}(w_{j})$ (Theorem \ref{algdirect}). By using a Darboux transformation between the weights $w_i$ and $w_j$, we explicitly obtain  a  generator of this module, i.e., we find a 
 differential operator $\mu_{i,j}$ such that 
$$\mathcal D(w_i,w_j)=\mu_{i,j}\cdot \mathcal{D}(w_{j}) .$$  
In particular, we prove in Theorem \ref{Dw_1w_2zero},  that $\mathcal D(w_i,w_j)=0$ if and only if  $w_i$ is not a Darboux transformation of  $w_j$.

  We use these results to determine when two direct sums of classical weights are Darboux transformations of each other. 
 In Theorem \ref{Darboux-direct sum}, we establish that  
 $$W(x)= w_1(x)\oplus  \cdots \oplus w_N(x)  \text{ is a  
 Darboux transformation of  } 
 \widetilde W(x)= \widetilde w_1(x)\oplus  \cdots \oplus \widetilde w_N(x),$$
 if and only if, up to permutations, each scalar weight $w_i$ is a Darboux transformation of $\widetilde w_i$.

In Section \ref{Sect-sph}, we consider an example of an irreducible weight matrix $W$ of Jacobi type, arising from matrix-valued spherical functions associated with the pair $(G, K)$ with $G=\mathrm{SU}(m+1)$ and $ K=\mathrm{U}(m)$.
We prove that this weight can be obtained as a Darboux transformation of a direct sum of classical Jacobi weights of the form $w_{\alpha+1,\beta}(x)\oplus w_{\alpha+1, \beta+1}(x)$ (then it cannot be a Darboux transformation of a single scalar weight $w_{\alpha,\beta}(x)I$). 
We obtain new results about the structure of its algebra $\mathcal D(W)$,  we prove that it is a commutative full algebra generated by the identity and two differential operators of order two.  

Finally, in Section \ref{7}, we present three examples of irreducible weight matrices obtained as Darboux transformations through the factorization of a certain differential operator in the algebra of a direct sum of classical weights.

\section{Matrix valued orthogonal polynomials and the algebra $\mathcal D(W)$}\label{sec-MOP}

 Let $W=W(x)$ be a weight matrix of size $N$ on the real line, that is, a complex $N\times N$ matrix-valued smooth function on a (possibly unbounded) interval $\mathcal{I}=(x_0,x_1)$ such that $W(x)$ is positive definite almost everywhere and with finite moments of all orders. We denote $\mathbb{W}(N,\mathcal{I})$  the set of all these weight matrices.

 Let $\operatorname{Mat}_N(\mathbb{C})$ be the algebra of all $N\times N$ complex matrices and let $\operatorname{Mat}_N(\mathbb{C})[x]$ be the algebra of polynomials in the indeterminate $x$ with coefficients in $\operatorname{Mat}_N(\mathbb{C})$. We consider the following Hermitian sesquilinear form in the linear space $\operatorname{Mat}_N(\mathbb{C})[x]$
\begin{equation*}
  \langle P,Q \rangle =  \langle P,Q \rangle_W = \int_{x_0}^{x_1} P(x) W(x) Q(x)^*\,dx.
\end{equation*}

Given a weight matrix $W$ one can construct sequences 
$\{Q_n\}_{n\in\mathbb{N}_0}$ of matrix-valued orthogonal polynomials, i.e., the $Q_n$ are polynomials of degree $n$ with nonsingular leading coefficient and $\langle Q_n,Q_m\rangle=0$ for $n\neq m$.
We observe that there exists a unique sequence of monic orthogonal polynomials $\{P_n\}_{n\in\mathbb{N}_0}$ in $\operatorname{Mat}_N(\mathbb{C})[x]$.
By following a standard argument (see \cite{K49} or \cite{K71}) one shows that the monic orthogonal polynomials $\{P_n\}_{n\in\mathbb{N}_0}$ satisfy a three-term recursion relation
\begin{equation}\label{ttrr}
    x P_n(x)=P_{n+1}(x) + B_{n}P_{n}(x)+ C_nP_{n-1}(x), \qquad n\in\mathbb{N}_0,
\end{equation}
where $P_{-1}=0$ and $B_n, C_n$ are matrices depending on $n$ and not on $x$.

Along this paper, we consider that an arbitrary matrix differential operator
\begin{equation}\label{D2}
  {D}=\sum_{i=0}^s \partial ^i F_i(x),\qquad \partial=\frac{d}{dx},
\end{equation}
acts on the right on a matrix-valued function $P $, i.e., $P(x)\cdot D=\sum_{i=0}^s \partial ^i (P)(x)F_i(x).$

We consider the algebra of these operators with polynomial coefficients $$\operatorname{Mat}_{N}(\Omega[x])=\Big\{D = \sum_{j=0}^{n} \partial^{j}F_{j}(x) \, : F_{j} \in \operatorname{Mat}_{N}(\mathbb{C})[x] \Big \}.$$
\noindent 
More generally, when necessary, we will also consider $\operatorname{Mat}_{N}(\Omega[[x]])$, the set of all differential operators with coefficients in $\mathbb{C}[[x]]$, the ring of power series with coefficients in $\mathbb{C}$.

Given a weight matrix $W$ and $\{P_n\}_{n\in \mathbb{N}_0}$  a sequence of matrix-valued orthogonal polynomials with respect to $W$, we introduce 
the algebra
of all differential operators that have the sequence $P_n$ as eigenfunctions, i.e., 
\begin{equation}\label{algDW}
  \mathcal D(W)=\left\{D\in \operatorname{Mat}_{N}(\Omega[x])\, : \, P_n \cdot D=\Lambda_n(D) P_n, \, \Lambda_n(D)\in \operatorname{Mat}_N(\mathbb{C}), \text{ for all }n\in\mathbb{N}_0\right\}
\end{equation}

We observe that the definition of $\mathcal D(W)$ depends only on the weight matrix $W$ and not  on the particular sequence of orthogonal polynomials,
since two sequences $\{P_n\}_{n\in\mathbb{N}_0}$ and $\{Q_n\}_{n\in\mathbb{N}_0}$ of matrix orthogonal polynomials are related by
$P_n=M_nQ_n$, with $\{M_n\}_{n\in\mathbb{N}_0}$ invertible matrices (see \cite{GT07}, Corollary 2.5).

The {\em formal adjoint} on $\operatorname{Mat}_{N}(\Omega([[x]])$, denoted  by $\mbox{}^*$, is the unique involution extending Hermitian conjugate on $\operatorname{Mat}_N(\mathbb{C})[x]$ and sending $\partial I$ to $-\partial I$. 
The {\em formal $W$-adjoint} of $ \mathfrak{D}\in \operatorname{Mat}_{N}(\Omega([x])$, or  the formal adjoint of $\mathfrak D$ with respect to $W(x)$  is the differential operator $\mathfrak{D}^{\dagger} \in \operatorname{Mat}_{N}(\Omega[[x]])$ defined
by
$$\mathfrak{D}^{\dagger}:= W(x)\mathfrak{D}^{\ast}W(x)^{-1},$$
where  $\mathfrak{D}^{\ast}$ is the formal adjoint of 
$\mathfrak D$. 
An operator $\mathfrak{D}\in \operatorname{Mat}_{N}(\Omega[x])$ is called {\em $W$-adjointable} if there exists 
$\widetilde {\mathfrak{D}} \in \operatorname{Mat}_{N}(\Omega[x])$, such that
$\langle P \cdot \mathfrak{D},Q\rangle=\langle P,Q \cdot \widetilde{\mathfrak{D}}\rangle,$ for all $P,Q\in \operatorname{Mat}_N(\mathbb{C})[x]$. Then we say that the operator $\widetilde{\mathfrak D}$ is the $W$-adjoint of $\mathfrak D $.

\smallskip
For the algebra $\mathcal{D}(W)$ we have 
the following result.

\begin{prop} \label{adjunta D(W)} 
    If $D \in \mathcal{D}(W)$, then  $D$ is $W$-adjointable. Moreover, there exists a unique $\tilde D$ in $\mathcal D(W)$ such that 
  $$ \langle P D, Q\rangle  = \langle P, Q\tilde D\rangle \quad \text{ for all } P, Q \in \operatorname{Mat}_{N}(\mathbb{C})[x]. $$
\end{prop}
\begin{proof}
    See Corollary 4.5 in \cite{GT07}.
\end{proof}

\begin{remark}
    To avoid pathological cases, we will consider only weight matrices $W(x)$ such that for each integer $n \geq 0$ the $n$-th derivative $W^{(n)}(x)$ decreases exponentially at infinity and there exists a scalar polynomial $p_{n}(x)$ such that $W^{(n)}(x)p_{n}(x)$ has finite moments. See \cite[Section 2.2]{CY18}.
\end{remark}

From \cite{CY18} we also obtain the following result:

\begin{prop}\label{adjuntas} 
    If $D \in \mathcal{D}(W)$, then $D^{\dagger}$ is the $W$-adjoint of $D$. 
\end{prop}

For $\mathfrak D= \sum_{j=0}^n \partial ^j F_j \in \operatorname{Mat}_{N}(\Omega[x]) $,  the formal $W$-adjoint of $\mathfrak D$ is given by $\mathfrak{D}^\dagger= \sum_{k=0}^n \partial ^k G_k$, with
\begin{equation}\label{daga}
    G_k=  \sum_{j = 0}^{n-k}(-1)^{n-j} \binom{n-j}{k}
    (WF_{n-j}^*)^{(n-k-j)} W^{-1} , \qquad \text{for } 0\leq k\leq n.  
 \end{equation}
 It is a matter of careful integration by parts to see that $\langle \, P \cdot \mathfrak{D}, Q\, \rangle=\langle \, P,\, Q \cdot {\mathfrak{D}}^\dagger \rangle$ if  the following set of ``boundary conditions"
  \begin{equation} 
  \lim_{x\to x_i}\, \sum_{j = 0}^{p-1}(-1)^{n-j+p-1} \binom{n-j}{k}
    \big (F_{n-j}(x)W(x)\big)^{(p-1-j)}=0 
      \end{equation}
are satisfied for $1\leq p\leq n$ and $0\leq k\leq n-p$, where $x_i$ are the endpoints of the support of the weight $W$. 
  
\ 

We say that a differential operator $D\in \mathcal D(W)$ is $W$-{\em symmetric} if $\langle P \cdot {D},Q\rangle=\langle P,Q \cdot D\rangle$, for all $P,Q\in \operatorname{Mat}_N(\mathbb{C})[x]$.
An operator $\mathfrak{D}\in \operatorname{Mat}_{N}(\Omega([x])$ is called {\em formally} $W${\em -symmetric} if $\mathfrak{D}^{\dagger} = \mathfrak{D}$.
In particular if $\mathfrak{D}\in \mathcal{D}(W)$, then $\mathfrak{D}$ is $W$-symmetric if and only if it is formally $W$-symmetric.

It is shown that the set $\mathcal S(W)$ of all symmetric operators in $\mathcal D(W)$ is a real form of the space $\mathcal D(W)$, i.e.,
\begin{equation}\label{symop}
    \mathcal D(W)= \mathcal S (W)\oplus i \,\mathcal S (W),
\end{equation}
as real vector spaces. 

The condition of symmetry for a 
differential operator in the algebra $\mathcal D(W)$ is equivalent to the following set of differential equations involving the weight $W$ and the coefficients of  $D$.

\begin{thm}\label{equivDsymm}
 Let $\mathfrak{D} =\sum_{i=0}^n \partial^i F_i(x)$ be a  differential operator of order $n$ in $\mathcal{D}(W)$.
 Then $\mathfrak{D}$ is $W$-symmetric if and only if
\begin{equation*}
    \sum_{j = 0}^{n-k}(-1)^{n-j} \binom{n-j}{k}
    (F_{n-j}W)^{(n-k-j)}= WF_{k}^\ast   
\end{equation*}
for all $0\leq k \leq n$.
 \end{thm}
 \begin{proof}  
 An operator $\mathfrak{D}\in \mathcal{D}(W)$ is $W$-symmetric if and only if $\mathfrak{D}= \mathfrak{D}^\dagger$. The statement follows by using the explicit expression of the coefficients of  $\mathfrak{D}^\dagger$ given in  \eqref{daga}. 
 \end{proof}

\begin{definition}
    Let $\mathcal{V} \in \operatorname{Mat}_{N}(\Omega [x])$, we call $\mathcal{V}$ \textbf{degree-preserving} if the degree of $P(x) \cdot \mathcal{V}$ is equal to the degree of $P(x)$ for all $P(x) \in \operatorname{Mat}_{N}(\mathbb{C})[x]$.
\end{definition}

Observe that a degree-preserving differential operator $\mathcal{V} = \sum_{j=0}^{m}\partial^{j}F_{j} \in \operatorname{Mat}_{N}(\Omega[x])$ satisfies $\deg(F_{j})\leq j$ for all $j$. 

Furthermore, it should be noted that every degree-preserving differential operator $\mathcal V$ is not a zero divisor in the algebra of differential operators $\operatorname{Mat}_{N}(\Omega[x])$, i.e., if $\mathcal V D=0 $ or $D\mathcal V=0$ then $D=0$. 

\begin{prop}\label{preserving nonsingular}
    If $\mathcal{V} \in \operatorname{Mat}_{N}(\Omega[x])$ is a degree-preserving  operator and 
    $P$ is a matrix polynomial with 
    %\(P(x) \in \operatorname{Mat}_{N}(\mathbb{C})[x]\) has 
    nonsingular leading coefficient,  then \(P(x) \cdot \mathcal{V}\) has a nonsingular leading coefficient.
\end{prop}

\begin{proof}
    Let \(A\) be the leading coefficient of \(P(x) \cdot \mathcal{V}\). If \(A\) is a singular matrix, then there exists a nonzero matrix \(B\) such that \(BA = 0\). Therefore, we get that \(BP(x) \cdot \mathcal{V}\) is of a degree less than that of \(BP(x)\), which is a contradiction because \(\mathcal{V}\) is a degree-preserving operator.
\end{proof}

\section{Darboux transformations}\label{Darboux-sect}

Throughout this section, we consider weight matrices $W$ and $\widetilde{W}$ supported on the same interval. We denote $P_n=P_{n}(x)$ and $\widetilde P_n=\widetilde{P}_{n}(x)$ their associated sequence of monic orthogonal polynomials, respectively.

\begin{definition}\label{darb def}
    We say that $\widetilde{W}$ is a Darboux transformation of $W$ if there exists a differential operator $D \in \mathcal{D}(W)$ that can be factorized as $D = \mathcal{V}\mathcal{N}$ with degree-preserving operators $\mathcal{V}, \, \mathcal{N} \in \operatorname{Mat}_{N}(\Omega[x])$, such that
    \begin{equation*}
            P_{n}(x) \cdot \mathcal{V}  = A_{n}\widetilde{P}_{n}(x) , \qquad  \text{ for all } n \in \mathbb N_0, 
    \end{equation*}
for some matrices $A_{n} \in \operatorname{Mat}_{N}(\mathbb{C})$.
\end{definition}

\begin{remark} 
In our definition of Darboux transformation, we require that the differential operators \(\mathcal{V}\) and \(\mathcal{N}\) be degree-preserving to ensure that \(P_n(x) \cdot \mathcal{V}\) forms a sequence of orthogonal polynomials with respect to \(\widetilde{W}\). Equivalently, we can require that the matrices \(A_n\) be nonsingular for all \(n \geq 0\).

While this imposes a greater restriction than usual in other contexts, such as considering exceptional orthogonal polynomials, our primary focus lies on weights whose sequences of orthogonal polynomials are eigenfunctions of differential operators.

\end{remark}

%Note that in the above definition, requiring \(\mathcal{V}\) to be degree-preserving is equivalent to require that the matrices $A_n$ be nonsingular, for all $n\geq 0$ and it  is equivalent to \(P_{n}(x) \cdot \mathcal{V}\)  be a sequence of orthogonal polynomials for \(\widetilde{W}\). 

%In fact, 
%This is because, by Proposition \ref{preserving nonsingular}, the \(A_{n}\) of the definition are nonsingular for all \(n \geq 0\). On the other hand, if \(\mathcal{V} \in \operatorname{Mat}_{N}(\Omega[x])\) is such that \(P_{n}(x) \cdot \mathcal{V} = A_{n} \widetilde{P}_{n}(x)\) with \(A_{n}\) nonsingular for all \(n \geq 0\), then \(\mathcal{V}\) is degree-preserving.

%\smallskip
In the next proposition, we see that the Darboux transformation is a symmetric property.

\begin{prop}\label{cons}
    Let $\widetilde{W}$ be a Darboux transformation of $W$ with $D = \mathcal{V}\mathcal{N}$ as in Definition \ref{darb def}. Then the differential operator $\widetilde{D} = \mathcal{N}\mathcal{V}$ belongs to $\mathcal{D}(\widetilde{W})$ and $W$ is a Darboux transformation of $\widetilde{W}$.
\end{prop}
\begin{proof}
    Since  $D=\mathcal V\mathcal N \in \mathcal D (W)$, we  have that  $P_n\cdot D=\Lambda_n(D)P_n$. Thus, for 
$\widetilde{D} = \mathcal{N}\mathcal{V}$ we have
\begin{equation*}
        \widetilde{P}_{n}(x)\cdot \widetilde{D} = A_{n}^{-1}P_{n}(x)\cdot \mathcal{V}\widetilde{D} = A_{n}^{-1}P_{n}(x)\cdot D\mathcal{V} = A_{n}^{-1}\Lambda_{n}(D)A_{n}\widetilde{P}_{n}(x).
\end{equation*}
Hence $\widetilde{D} \in \mathcal{D}(\widetilde{W})$, and 
\begin{equation}\label{1}
    \widetilde{P}_{n}(x)\cdot \mathcal{N} = A_{n}^{-1}\Lambda_{n}(D)P_{n}(x).
\end{equation}
Therefore $W$ is a Darboux transformation of $\widetilde{W}$.
\end{proof}

The Darboux transformations of weights are closely related to the algebra of differential operators.

\begin{prop}\label{art}
    Let $\widetilde{W}$ be a Darboux transformation of $W$, and let $\mathcal{V}$ and $\mathcal{N}$ be the operators given in Definition \ref{darb def}. Then the differential operator $\left(\begin{smallmatrix} 0 && \mathcal{V} \\ \mathcal{N} && 0 \end{smallmatrix}\right)$ belongs to $ \mathcal{D}(W \oplus \widetilde{W})$. 
\end{prop}
\begin{proof}
    Let $P_{n}(x)$ and $\widetilde{P}_{n}(x)$ be the sequences of monic orthogonal polynomials associated to $W$ and $\widetilde{W}$ respectively. 
    From Definition \ref{darb def} and  Proposition \ref{cons} we have that 
    $$P_{n}(x)\cdot \mathcal{V} = A_{n}\widetilde{P}_{n}(x) \quad \text{ and } \quad \widetilde{P}_{n}(x) \cdot \mathcal{N} = B_{n}P_{n}(x)$$
    for some invertible  matrices $A_{n}, \, B_{n} \in \operatorname{Mat}_{N}(\mathbb{C}).$ Thus $\left(\begin{smallmatrix} P_{n}(x) && 0 \\ 0 && \widetilde{P}_{n}(x) \end{smallmatrix}\right)$, the sequence of monic orthogonal polynomials associated to $W \oplus \widetilde{W}$, satisfies    
    $$\left( \begin{matrix} P_{n}(x) && 0 \\ 0 && \widetilde{P}_{n}(x) \end{matrix} \right) \cdot \left(\begin{matrix} 0 && \mathcal{V} \\ \mathcal{N} && 0 \end{matrix}\right) = \left(\begin{matrix} 0 && A_{n} \\ B_{n} && 0 \end{matrix} \right)\left(\begin{matrix} P_{n}(x) && 0 \\ 0 && \widetilde{P}_{n}(x) \end{matrix}\right).$$
    Therefore $\left(\begin{smallmatrix} 0 && \mathcal{V} \\ \mathcal{N} && 0 \end{smallmatrix}\right) $ belongs to $ \mathcal{D}(W \oplus \widetilde{W})$.
\end{proof}

\smallskip
Given two weight matrices $W$ and $\widetilde{W}$  supported on the same interval we introduce the relation 
\begin{equation}\label{equivrelat}
W \sim \widetilde{W}  \qquad \text{ if and only if }  \qquad \widetilde{W}  \text{ is a Darboux transformation of $W$}. 
\end{equation}

\smallskip
\begin{prop}\label{equivalence relation}
  The Darboux transformation defines an equivalence relation in the set of weight matrices in $\mathbb{W}(N,\mathcal{I})$.
\end{prop}
\begin{proof} 
     It is clear that  $W\sim W$ because the identity operator $I \in \mathcal{D}(W)$ can be factorized as $I = I\cdot I$. The symmetry follows directly from  Proposition \ref{cons}.

If $W_{1} \sim W_{2}$ and $W_{2} \sim W_{3}$, then there exist degree-preserving differential operators $\mathcal{N}_{1}$, $\mathcal{N}_{2}$, $\mathcal{V}_{1}$ and $\mathcal{V}_{2}$ such that $D_{1} = \mathcal{V}_{1}\mathcal{N}_{1} \in \mathcal{D}(W_{1})$, $D_{2} = \mathcal{V}_{2}\mathcal{N}_{2} \in \mathcal{D}(W_{2})$ with $P^{W_{1}}_{n}(x)\cdot \mathcal{V}_{1} = A_{n}P^{W_{2}}_{n}(x)$ and $P^{W_{2}}_{n}(x) \cdot \mathcal{V}_{2} = B_{n} P^{W_{3}}_{n}(x)$ for some invertible matrices $A_{n}$, $B_{n}$. Then the degree-preserving operators $\mathcal{V} = \mathcal{V}_{1}\mathcal{V}_{2}$ and $\mathcal{N} = \mathcal{N}_{2}\mathcal{N}_{1}$ satisfy 
    $$P^{W_{1}}_{n}(x) \cdot \mathcal{V} = A_{n}P^{W_{2}}_{n}(x)\cdot \mathcal{V}_{2} = A_{n}B_{n}P^{W_{3}}_{n}(x).$$
It remains to see that $D = \mathcal{V} \mathcal{N} = \mathcal{V}_{1}D_{2}\mathcal{N}_{1} \in \mathcal{D}(W_{1})$.  In fact, 
    $$P^{W_{1}}_{n}(x) \cdot D = A_{n}P^{W_{2}}_{n}(x)\cdot D_{2}\mathcal{N}_{1} = A_{n}\Lambda_{n}(D_{2})P^{W_{2}}_{n}(x)\cdot \mathcal{N}_{1}$$ 
    and $P^{W_{2}}_{n}(x)\cdot \mathcal{N}_{1} = A_{n}^{-1}\Lambda_{n}(D_{1})P^{W_{1}}_{n}(x)$. Therefore  $W_{1} \sim W_{3}$.
    \end{proof}
    
\begin{definition}  
    We say that $\widetilde{W}$ is Darboux-equivalent to $W$ if  $\widetilde{W}$ is a Darboux transformation of $W$.
    We denote by $[W]$ the Darboux-equivalence class of $W$.
\end{definition}

\begin{prop} 
Equivalent weights are Darboux-equivalent weights.
\end{prop}
\begin{proof} 
    If $W$ and  $\widetilde{W}$ are equivalent weights, 
    then there exists an invertible matrix $M \in \operatorname{Mat}_{N}(\mathbb{C})$ such that $$\widetilde{W}(x) = MW(x)M^{\ast},\quad \text{ for all $x$. }$$  
    Thus we have that $\widetilde{P}_{n}(x) M = MP_{n}(x)$ for all $n\geq 0$. Then we can factorize the zero-order differential operator $I = M M^{-1}\in \mathcal{D}(W)$. 
    Hence $W$ is a Darboux transformation of $\widetilde{W}$.
\end{proof}

\begin{prop} \label{direct}
If $\widetilde{W}_{1}$ and $\widetilde{W}_{2}$ are Darboux transformations of $W_{1}$ and $W_{2}$, respectively, then $\widetilde{W}_{1} \oplus \widetilde{W}_{2}$ is a Darboux transformation of $W_{1} \oplus W_{2}$.
\end{prop}
\begin{proof}
    For $j=1,2$ let $P^{W_{j}}_{n}(x)$, $\widetilde{P}^{W_{j}}_{n}(x)$ be the sequence of monic orthogonal polynomials for $W_{j}$ and $\widetilde{W}_{j}$, respectively. 
 
 There exist degree-preserving differential operators $\mathcal{V}_{j}$ and $\mathcal{N}_{j}$ such that  $D_{j} = \mathcal{V}_{j} \mathcal{N}_{j} \in \mathcal{D}(W_{j})$, $P^{W_{1}}_{n}(x) \cdot \mathcal{V}_{1} = A_{n}\widetilde{P}^{W_{1}}_{n}(x)$, and $P^{W_{2}}_{n}(x) \cdot \mathcal{V}_{2} = B_{n}\widetilde{P}^{W_{2}}_{n}(x)$ for some  nonsingular matrices $A_{n} \in \operatorname{Mat}_{N_{1}}(\mathbb{C})$, $B_{n} \in \operatorname{Mat}_{N_{2}}(\mathbb{C})$.  
 The operators $\mathcal{V} = \operatorname{diag}(\mathcal{V}_{1},\mathcal{V}_{2})$ and  $\mathcal{N} = \operatorname{diag}(\mathcal{N}_{1}, \mathcal{N}_{2})$ are  degree-preserving differential operators, $D =\mathcal{V}\mathcal{N} = \operatorname{diag}(D_{1},D_{2}) \in \mathcal{D}(W_{1} \oplus W_{2})$ and 
 $$\operatorname{diag}(P^{W_{1}}_{n}(x),P^{W_{2}}_{n}(x))\cdot \mathcal{V} = \operatorname{diag}(A_{n},B_{n}) \operatorname{diag}(\widetilde{P}^{W_{1}}_{n}(x),\widetilde{P}^{W_{2}}_{n}(x)).$$  
 Therefore $\widetilde{W}_{1} \oplus \widetilde{W}_{2}$ is a Darboux transformation of $W_{1} \oplus W_{2}$.
\end{proof}

The algebra of differential operators $\mathcal{D}(W)$ is closely related to the algebra $\mathcal{D}(\widetilde W)$ when $W$ and $\widetilde W$ are Darboux-equivalent weights.

\begin{prop} \label{alg debil}
    Let 
    $\widetilde{W}$ be a Darboux transformation of $W$ and let   $D = \mathcal{V}\mathcal{N} \in \mathcal{D}(W)$  as in Definition \ref{darb def}. Then  
    $$\mathcal{V}\mathcal{D}(\widetilde{W})\mathcal{N} \subseteq \mathcal{D}(W) \quad \text{ and } \quad \mathcal{N} \mathcal{D}(W)\mathcal{V} \subseteq \mathcal{D}(\widetilde{W}).$$
\end{prop}
\begin{proof}
    Let $P_{n}(x)$ and $\widetilde{P}_{n}(x)$ be the sequence of monic orthogonal polynomials asociated to $W$ and $\widetilde{W}$ respectively. From  Proposition \ref{cons} we have that $P_{n}(x)\cdot \mathcal{V} = A_{n} \widetilde{P}_{n}(x)$ and $\widetilde{P}_{n}(x) \cdot \mathcal{N} = B_{n}P_{n}(x)$, for some invertible matrices $A_{n}$ and $B_{n}$. 
    
For  $\mathcal{A} \in \mathcal{D}(\widetilde W)$ and  $n\geq 0$, it is enough to prove that 
$P_{n}(x)$ is an eigenfunction of $\mathcal{V}\mathcal{A}\mathcal{N}$.  
 In fact    $$P_{n}(x) \cdot \mathcal{V}\mathcal{A}\mathcal{N} = A_{n}\widetilde{P}_{n}(x)\cdot \mathcal{A}\mathcal{N} =A_{n} \Lambda_{n}(\mathcal{A})\widetilde{P}_{n}(x)\cdot \mathcal{N} =  A_{n} \Lambda_{n}(\mathcal{A})B_{n}P_{n}(x).$$
    Therefore  $\mathcal{V}\mathcal{A}\mathcal{N}$ belongs to $\mathcal{D}(W)$. In the same way we obtain that $\mathcal{N} \mathcal{D}(W)\mathcal{V} \subseteq \mathcal{D}(\widetilde{W})$. 
\end{proof}

We now obtain some important consequences of the previous results.

\begin{thm}\label{conmutativa}
If $\widetilde{W}$ is a Darboux transformation of $W$ then $\mathcal{D}(W)$ is a commutative algebra if and only if $\mathcal{D}(\widetilde{W})$ is a commutative algebra.
\end{thm}
\begin{proof}
  Let $\mathcal V$ and $\mathcal N$ be degree-preserving operators as in Definition \ref{darb def}. We have $D = \mathcal{V}\mathcal{N}\in \mathcal{D}(W)$   and  $\widetilde{D} = \mathcal{N}\mathcal{V}\in \mathcal{D}(\widetilde{W})$. From Proposition \ref{alg debil} we have that  $\mathcal{V}\mathcal{A}\mathcal{N} \in \mathcal{D}(W)$ for all $\mathcal{A} \in \mathcal{D}(\widetilde{W})$.

We observe that if  $\mathcal{D}(W)$ is a commutative algebra then  $\widetilde{D}$ belongs to the center of $\mathcal{D}(\widetilde{W})$. In fact, for $\mathcal{A} \in \mathcal{D}(\widetilde{W})$, we get  
$$\mathcal{V}\mathcal{A}\widetilde{D}\mathcal{N} = \mathcal{V}\mathcal{A}\mathcal{N}D = D\mathcal{V}\mathcal{A}\mathcal{N} = \mathcal{V}\widetilde{D}\mathcal{A}\mathcal{N},$$
since $\mathcal{V}$ and $\mathcal{N}$ are degree-preserving operators they are not zero divisors, thus $\mathcal{A}\widetilde{D} = \widetilde{D}\mathcal{A}$.
    
    Now, for $\mathcal{A}_{1}$, $\mathcal{A}_{2} \in \mathcal{D}(\widetilde{W})$ we have that $(\mathcal{V}\mathcal{A}_{1}\mathcal{N})(\mathcal{V}\mathcal{A}_{2}\mathcal{N}) = (\mathcal{V}\mathcal{A}_{2}\mathcal{N})(\mathcal{V}\mathcal{A}_{1}\mathcal{N}),$
thus  $$\mathcal{V}\mathcal{A}_{1}\widetilde{D}\mathcal{A}_{2} \mathcal{N} = \mathcal{V}\mathcal{A}_{2}\widetilde{D}\mathcal{A}_{1}\mathcal{N}$$
and by using that $\widetilde{D}$ belongs to the center of $\mathcal{D}(\widetilde{W})$ and $\mathcal{V}$ and $\mathcal{N}$ are not zero divisors, we obtain that $\mathcal{A}_{1}\mathcal{A}_{2} = \mathcal{A}_{2}\mathcal{A}_{1}$.
\end{proof}

Later, in Section \ref{sub}, we will examine the algebras $\mathcal{D}(W)$ for a direct sum of scalar weights. Specifically, in Proposition \ref{com} we will characterize the conditions under which this algebra is commutative.

\begin{thm}\label{com 2}
    If a weight matrix $W$ is a Darboux transformation of 
    a direct sum of scalar weights, 
    $%\widetilde W=
    w_{1} \oplus  \cdots \oplus w_{N}$, then the algebra $\mathcal{D}({W})$ is commutative if and only if $w_{i}$ is not a Darboux transformation of $w_{j}$ for all $i \not= j$.
\end{thm}

\section{Classical scalar weights}\label{classical-section}
In the real line $\mathbb R$, we know that the only weights $w$  containing a second-order operator $\delta$ in the algebra $\mathcal D(w)$ are, 
up to an affine change of coordinates, the classical  Hermite, Laguerre, and Jacobi weights. The differential operator $\delta $ is the classical Hermite, Laguerre, or Jacobi operator, respectively.
 We refer to them as the {\em classical scalar weights}. 

\begin{table}[H]
\begin{center} \small
    \begin{tabular}{|c|c|c|c|}
     \hline    Type & differential operator  & weight   & support    \\    \hline &&& \\ shifted Hermite   &   $\delta =  \partial^{2} -2 \partial (x-b)$ & $w_b(x)=e^{-x^2+2bx}$ & $(-\infty, \infty)$   \\ \hline & & & \\   Laguerre    &  $\delta =  \partial^{2}x  + \partial(\alpha+1-x)$   & $w_\alpha(x)= e^{-x} x^\alpha$ & $ (0,\infty )$ \\ \hline  & & & \\
     Jacobi  & $\delta = \partial^{2} (1-x^2)  + \partial (\beta-\alpha-x(\alpha+\beta+2))$ & $w_{\alpha,\beta}(x) = (1-x)^\alpha (1+x)^\beta $ & $(-1,1)$ \\ & & & \\\hline
    \end{tabular}\caption{Classical scalar weights. } \label{weights-table}
    \label{weights-oper}
\end{center} 
\end{table}

\begin{remark} 
    Throughout this paper, we will consider the shifted Hermite weights because they are relevant for the construction of several weight matrices which are solutions to the Bochner problem and they are not Darboux transformations of classical weights. See \cite{BP23} and \cite{BP24-1}. 
\end{remark}

\begin{thm}(See \cite{M05}) 
Let $w(x)$ be a classical weight (Hermite, Laguerre, or Jacobi). Then the algebra $\mathcal D(w) $
is a polynomial algebra in the differential operator $\delta$, i.e., $\mathcal D(w)= \mathbb C[\delta]$.
\end{thm}

In all these classical cases, the structure of the algebra $\mathcal D(w)$ is very well known. 

We want to study the Darboux-equivalence classes between the scalar classical weights of the same type. First, we will present some general results that can be derived by studying the algebra of the direct sum of scalar weights.

\begin{prop} \label{diag D}
    Let $w_{1}$ and $w_{2}$ be scalar weights such that $w_{1}(x)w_{2}(x)^{-1}$ is not a rational function.  Then the algebra of the direct sum $\mathcal{D}(w_{1} \oplus w_{2})$ is a diagonal algebra, that is $$\mathcal{D}(w_{1}\oplus w_{2}) \subseteq \left \{ \sum_{j=0}^{n} \partial^{j} F_{j}(x) \, :
\, 
F_{j}(x) = \begin{pmatrix} p_{j}(x) && 0 \\ 0 && q_{j}(x) \end{pmatrix}, p_{j},q_{j} \in \mathbb{C}[x], \, n\geq 0 \right \}.$$
\end{prop}

\begin{proof}
Let $W(x) = w_{1}(x) \oplus w_{2}(x)$. 
From \eqref{symop} 
it is enough to prove that $\mathcal{S}(W)$ is a diagonal algebra.
Let  $D=\sum_{j=0}^{n} \partial^{j}F_{j}(x)$ be a $W$-symmetric operator in $\mathcal D(W)$.  
The coefficients $F_j(x)=\begin{pmatrix} p_{j}(x) && r_{j}(x) \\ g_{j}(x) && q_{j}(x) \end{pmatrix} \in \operatorname{Mat}_{2}(\mathbb{C})[x]$ satisfy the following set of equations given in Proposition \ref{equivDsymm} 
\begin{equation}\label{GG}
    F_k(x)=  \sum_{j = 0}^{n-k}(-1)^{n-j} \binom{n-j}{k}
    (W(x)F_{n-j}(x)^*)^{(n-k-j)}W(x)^{-1} , \qquad \text{for } 0\leq k\leq n.  
 \end{equation}
 In particular, for $k=n$ we get  
$$F_{n}(x) =(-1)^{n} W(x)F_{n}(x)^{\ast}W(x)^{-1} =(-1)^{n} \begin{pmatrix} \overline{p}_{n}(x) && w_{1}(x)w_{2}(x)^{-1}\overline{g}_{n}(x) \\ w_{2}(x)w_{1}(x)^{-1}\overline{r}_{n}(x) && \overline{q}_{n}(x)  \end{pmatrix},$$ 
which implies that 
$g_{n}(x)=0$ and $r_{n}(x)=0$ and hence  $F_n$ is a diagonal matrix. 
By induction on $k$ we can assume that $F_{n-j}$ is diagonal, for all $0 \leq j \leq k$. 
From \eqref{GG} we have that 
$$F_{n-k-1}  = (-1)^{n-k-1} F_{n-k-1}^*+\sum_{j=0}^{k}(-1)^{n-j}\binom{n-j}{n-k-1}(WF_{n-j}^{\ast})^{(k+1-j)}W^{-1}.$$
Since the weight $W$ is diagonal, we get 
$$\sum_{j=0}^{k}(-1)^{n-j}\binom{n-j}{n-k-1}(W(x)F_{n-j}(x)^{\ast})^{(k+1-j)}W(x)^{-1} = \begin{pmatrix} \alpha(x) && 0 \\ 0 && \beta(x) \end{pmatrix},$$
for some polynomials  $\alpha,\beta$.  
Thus 
\begin{equation*}
    \begin{split}
        F_{n-k-1} &  = \begin{pmatrix} \alpha(x) + (-1)^{n-k-1}\overline{p}_{n-k-1}(x) && (-1)^{n-k-1}w_{1}(x)w_{2}(x)^{-1}\overline{g}_{n-k-1}(x) \\ (-1)^{n-k-1}w_{2}(x)w_{1}(x)^{-1}\overline{r}_{n-k-1}(x) && \beta(x) + (-1)^{n-k-1}\overline{q}_{n-k-1}(x) \end{pmatrix}
    \end{split}
\end{equation*}
is a matrix polynomial. Thus $g_{n-k-1}(x)=0$ and $r_{n-k-1}(x)=0$ and this concludes the proof.
\end{proof}

\begin{prop}\label{clas dis}
   Let $w_{1}(x)$ and $w_{2}(x)$ be scalar weights such that $w_{1}(x)w_{2}(x)^{-1}$ is not a rational function. 
   Then the weight $w_{2}$ is not a Darboux transformation of $w_{1}$.
\end{prop}
\begin{proof}
    Let us assume that $w_{2}$ is a Darboux transformation of $w_{1}$ and let $\mathcal{V}$ and $\mathcal{N}$ be some degree-preserving differential operators as given in  Definition \ref{darb def}. Then by Proposition \ref{art}, the differential operator 
    $\left(\begin{smallmatrix} 0 && \mathcal{V} \\ \mathcal{N} && 0 \end{smallmatrix}\right) $ belongs to $\mathcal{D}(w_{1} \oplus w_{2})$, which is a diagonal algebra by Proposition \ref{diag D}.  Therefore $\mathcal{V}  = \mathcal{N}=0$, which contradicts the fact that $\mathcal{V}$ and $\mathcal{N}$ are degree-preserving.
\end{proof}

\section{Darboux-equivalence classes for scalar weights} \label{sub}

This section aims to identify the Darboux-equivalent classes of the classical Hermite, Laguerre, and Jacobi weights.

\smallskip
The Darboux-equivalence class of a (shifted) Hermite weight $w_{a}(x) = e^{-x^{2} + 2ax}$ consists only of the weight $w_a$. In fact, $w_{a}(x) w_{b}(x)^{-1}$ is not a rational function  when $a\neq b$ and thus, by Proposition \ref{clas dis}, they are not Darboux-equivalent.

\smallskip
The case of Laguerre weights $w_{\alpha}(x) = e^{-x}x^{\alpha}$, $\alpha > - 1$, is more interesting. As we will see in the next theorem, the set of all Darboux-equivalence classes is parameterized by the interval $[0,1)$.

\begin{thm} \label{laguerre scalar equivalences}
    Let $\alpha, \, \beta > -1$. 
    The Laguerre weight $w_{\alpha}$ is a Darboux transformation of $w_{\beta }$ if and only if $  \alpha -\beta \in \mathbb{Z}$.
\end{thm}
\begin{proof}
    If $\alpha -\beta \not \in \mathbb{Z}$ then $w_{\alpha}(x)w_{\beta}(x)^{-1} = x^{\alpha-\beta}$ is not a rational function and  thus $w_\alpha$ is not a Darboux transformation of $w_\beta$ by Proposition \ref{clas dis}.
    For the converse assertion, it is enough to prove that  $w_{\alpha+1}$ is a Darboux transformation of $w_{\alpha}$ for any   $\alpha > -1$. 
    The operator $\tilde{\delta}_{\alpha} = \delta_{\alpha} -\alpha - 1 \in \mathcal{D}(w_{\alpha})$, where $\delta_{\alpha} = \partial^{2}x + \partial(\alpha+1-x)$ is the classical second-order differential operator of Laguerre, can be factorized as $\tilde{\delta}_{\alpha} = \mathcal{V} \mathcal{N}$, with $\mathcal{V} = \partial - 1$ and $\mathcal{N} = \partial x + \alpha + 1$. 
    
    The monic Laguerre  polynomials $L_{n}^{\alpha}(x)$, associated to the weight $w_{\alpha}$, satisfy the identity 
    $$L_{n}^{\alpha}(x)\cdot \mathcal{V} = - L_{n}^{\alpha+1}(x), $$  
    (see, for example  (5.1.13) and (5.1.14) in \cite{Sz}, p.102).  
    Therefore the weight $w_{\alpha+1}$ is a Darboux transformation of $w_{\alpha}$.
\end{proof}

We now study the scalar Jacobi weights $w_{\alpha,\beta}(x) = (1-x)^{\alpha}(1+x)^{\beta}$, for $\alpha, \, \beta > - 1$. 
In Theorem \ref{Jacobi}, we will see that the scalar Jacobi weights $w_{\alpha,\beta}$ and $w_{\widetilde{\alpha},\widetilde{\beta}}$ are in the same Darboux-equivalence class if and only if their parameters satisfy $\widetilde{\alpha} = \alpha + k$ and $\widetilde{\beta} = \beta - k$ for some integer $k$. In particular, there is a bijective correspondence between the equivalence class representatives and the set $(-1,\infty)\times (-1,0]$.

\begin{prop}\label{pip}
    Let $\alpha, \, \beta, \, \widetilde{\alpha}, \, \widetilde{\beta} > -1$ such  that $\widetilde{\alpha}-\alpha \notin \mathbb{Z}$ or $\widetilde{\beta}-\beta \notin \mathbb{Z}$. Then $w_{\widetilde{\alpha},\widetilde{\beta}}$ is not a Darboux transformation of $w_{\alpha,\beta}$.
\end{prop}
\begin{proof}
    If $\widetilde{\alpha}-\alpha \notin \mathbb{Z}$ or $\widetilde{\beta}-\beta \notin \mathbb{Z}$, then $w_{\alpha,\beta}(x)w_{\widetilde{\alpha},\widetilde{\beta}}(x)^{-1}$ is not a rational function, therefore by Proposition \ref{clas dis} $w_{\alpha,\beta}$ is not Darboux-equivalent to $w_{\widetilde{\alpha},\widetilde{\beta}}$.
\end{proof}

\begin{prop}\label{pop}
  For   $m \in \mathbb{N}$, the Jacobi weight $w_{\alpha+m,\beta}$ is not a Darboux transformation of $w_{\alpha,\beta}$.
\end{prop}
\begin{proof}
   Any differential operator in $\mathcal{D}(w_{\alpha,\beta})$ is  of even order because  the algebra $\mathcal{D}(w_{\alpha,\beta})$ is a polynomial algebra in
   $\delta_{\alpha,\beta} = \partial^{2}(1-x^{2}) + \partial (\beta - \alpha - x(\alpha+\beta+2))$.
   
    By induction on $n$, we have that 
    $$\delta_{\alpha,\beta}^{n}= \partial^{2n}(1-x^{2})^{n} + \partial^{2n-1}(1-x^{2})^{n-1}n(\beta-\alpha-x(\alpha+\beta+2n))+\text{lower order operators}. $$  

    If  $w_{\alpha+m,\beta}$ is a Darboux transformation of $w_{\alpha,\beta}$, then there exist degree-preserving differential operators $\mathcal{V} = \sum_{i=0}^{s}\partial^{i}f_{i}(x)$, $\mathcal{N} = \sum_{j=0}^{r}\partial^{j}g_{j}(x)$ such that $ \mathcal{V}\mathcal{N} \in \mathcal{D}(w_{\alpha,\beta})$ and $\mathcal{N}\mathcal{V} \in \mathcal{D}(w_{\alpha+m,\beta})$. In particular, $f_{i}$ and $ g_{i}$ are polynomials of degree less than or equal to $i$ and we have 
    
  \begin{align*}
    \mathcal{V}\mathcal{N} &=\sum_{i=0}^{s}\partial^{i}f_{i}(x)\sum_{j=0}^{r}\partial^{j}g_{r}(x) =  \sum_{k=0}^{n}c_{k}\,\delta_{\alpha,\beta}^{k} 
    \\
    \mathcal{N}\mathcal{V} 
    &= \sum_{j=0}^{r}\partial^{j}g_{r}(x)\sum_{i=0}^{s}\partial^{i}f_{i}(x) = \sum_{\ell = 0}^{n} d_{\ell}\, \delta_{\alpha+m,\beta}^{\ell},
  \end{align*}
    for some $c_{k}, \, d_{\ell} \in \mathbb{C}$ and $n\geq 1$. 
   By comparing the two higher-order terms of the operators 
   $ \mathcal{V}\mathcal{N}$ and $ \mathcal{N}\mathcal{V}$ 
   we obtain that 
\begin{align}
    r + s  = 2n, \qquad c_{n}  &= d_{n}, \qquad f_{s}(x)g_{r}(x)  = c_{n}(1-x^{2})^{n}, \label{eq1}\\ 
    rg_{r}(x)f_{s}'(x) + g_{r}(x)f_{s-1}(x) + g_{r-1}(x)f_{s}(x) & =n c_{n} (1-x^{2})^{n-1}(\beta-\alpha-x(\alpha+\beta+2n)), \label{eq2}\\
    sf_{s}(x)g_{r}'(x) + f_{s}(x)g_{r-1}(x) + f_{s-1}(x)g_{r}(x) & = n c_{n} (1-x^{2})^{n-1}(\beta-\alpha-m-x(\alpha+m+\beta+2n)). \label{eq3}
\end{align}

    Without loss of generality, we can assume that $c_{n} = 1$. 
    Thus, from \eqref{eq1} we  get 
    $$g_{r}(x) = (1-x)^{r_{1}}(1+x)^{r_{2}} \quad \text { and }  \quad f_{s}(x) = (1-x)^{n-r_{1}}(1+x)^{n-r_{2}},$$ with 
    $r = r_{1}+r_{2}$,  $s = 2n-r$ and $0 \leq r_{1},r_{2} \leq n$. 
    By subtracting \eqref{eq2} and \eqref{eq3} we obtain
    \begin{equation}\label{eq4}
    rg_{r}(x)f_{s}'(x)-sf_{s}(x)g_{r}'(x) = nm (1-x^{2})^{n-1}(1+x).
    \end{equation}
    But the left-hand side of \eqref{eq4}  is equal to $2n(r_{1}-r_{2})(1-x^{2})^{n-1} $, which is a contradiction because the right-hand side is a polynomial of degree $2n-1$.
   \end{proof}

   \begin{thm} \label{Jacobi}
   The Jacobi weight $w_{\alpha,\beta}$ is a Darboux transformation of  $w_{\widetilde{\alpha},\widetilde{\beta}}$ 
   if and only if $$\alpha+\beta =\widetilde{\alpha} + \widetilde{\beta}\quad \text{ and } \quad \alpha-\widetilde\alpha \in \mathbb{Z}.$$ 
    \end{thm} 
    \begin{proof}
        We first prove that $w_{\alpha+1,\beta}$ is a Darboux transformation of $w_{\alpha,\beta+1}$. Let $J_{n}^{\alpha,\beta}(x)$ be the sequence of monic Jacobi polynomials, which are orthogonal with respect to $w_{\alpha,\beta}$. 
 The differential operator $\tilde{\delta}_{\alpha,\beta+1} = \delta_{\alpha,\beta+1} - (\alpha+1)(\beta+1)$ belongs to $\mathcal{D}(w_{\alpha,\beta+1})$, where $\delta_{\alpha,\beta+1} = \partial^{2} (1-x^{2}) + \partial (\beta +1 - \alpha -(\alpha+\beta+3)x)$ is the classical second-order differential operator of the Jacobi weight, and it can be factorized  as $\tilde{\delta}_{\alpha,\beta+1} = \mathcal{V}\mathcal{N}$ with 
 $$\mathcal{V} = \partial (1+x) + \beta + 1 \quad \text{ and } \quad  \mathcal{N} = \partial (1-x) - \alpha - 1.$$ 
 
\noindent
 %Jacobi polynomials can be writt in terms of the hypergeometric function, 
 We can use Gauss' contiguous relations to obtain that 
\begin{equation*}
    (1+x)  \frac{d}{dx} J_n^{\alpha,\beta+1}(x) = -(\beta+1) J_n^{\alpha,\beta+1}(x) +(n+\beta+1) J_n^{\alpha+1,\beta}(x),
\end{equation*}
hence it follows that $J_{n}^{\alpha,\beta+1}(x) \cdot \mathcal{V} = (n+\beta+1) J_{n}^{\alpha+1,\beta}(x)$. Therefore  $w_{\alpha+1,\beta}$ is a Darboux transformation of $w_{\alpha,\beta+1}$.

        If $\widetilde{\alpha} = \alpha - k$ and $\widetilde{\beta} = \beta + k$, with $k\in \mathbb N$, then $w_{\alpha,\beta} \sim w_{\alpha - 1 , \beta + 1}  \sim \cdots \sim w_{\alpha-k,\beta+k} = w_{\widetilde{\alpha},\widetilde{\beta}}.$
        Therefore $w_{\alpha,\beta} \sim w_{\widetilde{\alpha},\widetilde{\beta}}$.
        
        Conversely, suppose that $w_{\widetilde{\alpha},\widetilde{\beta}}$ is a Darboux transformation of $w_{\alpha,\beta}$. By Proposition \ref{pip} it follows that $\widetilde{\alpha}-\alpha$ and $\widetilde{\beta}-\beta$ are integer numbers. 
           If  $\alpha - \widetilde{\alpha} > \widetilde{\beta} - \beta$ then $\alpha = \widetilde{\alpha} + k + m$ and $ \widetilde{\beta}=\beta  + k$ for some $k \in \mathbb{Z}$ and $m\in \mathbb{N}$. Thus $ w_{\alpha,\beta} = w_{\widetilde{\alpha}+k+m,\widetilde{\beta}-k}$ is a Darboux transformation of the weight $w_{\widetilde{\alpha}+m,\widetilde{\beta}}$, 
        which contradicts  Proposition \ref{pop}. 
        For  $\alpha-\widetilde{\alpha} < \widetilde{\beta}-\beta$ we proceed in the same way. Thus $\alpha - \widetilde{\alpha} = \widetilde{\beta}- \beta \in \mathbb{Z}$. 
    \end{proof}

\begin{remark} 
    For $\alpha, \, \beta > -1$, the Darboux-equivalence class of a Jacobi weight $w_{\alpha,\beta}$ is $$[w_{\alpha,\beta}] = \{ w_{\alpha+k,\beta-k} \, | \, k \in \mathbb{Z}, \, -\alpha-1 < k < \beta + 1 \}.$$
\end{remark}

\section{The algebra of a direct sum of classical scalar weights}\label{algstruct-sect}
This section aims to study and describe the algebra $\mathcal D(W)$ when $W$ is a direct sum of classical scalar weights of the same type 
$$W = w_{1}\oplus w_{2} \oplus \cdots \oplus w_{N},$$ 
where all $w_j$, $1\leq j\leq N$, are in the same family of (shifted) Hermite, Laguerre, or Jacobi weights.

Following \cite{TZ18}, given two scalar weights $w_{1}$, $w_{2}$ and $p^{w_{1}}_{n}(x),$ $p^{w_{2}}_{n}(x)$ their respective sequence of monic orthogonal polynomials, we define the set
\begin{equation} \label{module-def}
    \mathcal{D}(w_{1},w_{2}) = \{ \mathcal{T} \in \Omega[x] \, : \, p^{w_{1}}_{n}(x) \cdot \mathcal{T} = \Lambda_{n}(\mathcal{T})p^{w_{2}}_{n}(x), \, \Lambda_{n}(\mathcal{T}) \in \mathbb{C}, \, n\geq 0 \}. 
\end{equation}

It follows that $\mathcal{D}(w_{1},w_{2})$ is a left $\mathcal{D}(w_{1})$-module and a right $\mathcal{D}(w_{2})$-module. Furthermore, we have $\mathcal{D}(w_{1},w_{2})\mathcal{D}(w_{2},w_{1}) \subseteq \mathcal{D}(w_{1})$.
In particular, if $w_{1} = w_{2}$, then $\mathcal{D}(w_{1},w_{2}) = \mathcal{D}(w_{1})$. 
Therefore
$$\mathcal{D}(w_{1} \oplus w_{2} ) = \begin{pmatrix} \mathcal{D}(w_{1}) && \mathcal{D}(w_{1},w_{2}) \\ \mathcal{D}(w_{2},w_{1}) && \mathcal{D}(w_{2}) \end{pmatrix}.$$

\begin{prop}\label{alg 22} 
    Let $w_{1}, \, w_{2} $ be scalar weights. Then  $\mathcal{D}(w_{1},w_{2}) \not= 0$ if and only if $\mathcal{D}(w_{2},w_{1})\not= 0$.
\end{prop}
\begin{proof}
    If  $0\neq \mathcal{T} \in \mathcal{D}(w_{1},w_{2})$, then $D = \left(\begin{smallmatrix} 0 && \mathcal{T} \\ 0 && 0 \end{smallmatrix} \right) \in \mathcal{D}(w_{1} \oplus w_{2})$. By Proposition \ref{adjuntas}, it follows that $D^{\dagger}$, the adjoint of $D$,  belongs to $\mathcal{D}(w_{1} \oplus w_{2})$ and is of the form
    $$D^{\dagger} =\left(\begin{smallmatrix} w_{1}(x) && 0 \\ 0 && w_{2}(x) \end{smallmatrix} \right) \left(\begin{smallmatrix} 0 && \mathcal{T} \\ 0 && 0 \end{smallmatrix} \right)^{\ast} \left(\begin{smallmatrix} w_{1}(x)^{-1} && 0 \\ 0 && w_{2}(x)^{-1} \end{smallmatrix} \right) = \left(\begin{smallmatrix} 0 && 0 \\ w_{1}(x)\mathcal{T}^{\ast}w_{2}(x)^{-1} && 0 \end{smallmatrix} \right)  \in \mathcal{D}(w_{1} \oplus w_{2}). $$
    Hence $ 0 \not = \widetilde{\mathcal{T}} = w_{1}(x)\mathcal{T}^{\ast}w_{2}(x)^{-1} \in \mathcal{D}(w_{2},w_{1})$. 
\end{proof}

Given   $w_{1}(x),\ldots,w_{N}(x)$ 
a collection of classical scalar weights of the same type (shifted Hermite, Laguerre, or Jacobi), we consider the diagonal weight matrix  $W= w_{1}\oplus \cdots \oplus w_{N}$. We have 
\begin{equation}\label{zu}
    \mathcal{D}(w_{1}\oplus \cdots \oplus w_{N} ) = \sum_{j=1}^{N}\sum_{k=1}^{N}\mathcal{D}(w_{j},w_{k})E_{j,k},
\end{equation}
where $E_{i,j}$ denotes the $N\times N$ matrix with the $(i,j)$-entry equal to $1$ and equal to $0$ elsewhere.

\smallskip

Observe that the $j$-th diagonal entry of a differential operator $D\in \mathcal D(W)$ belongs to the algebra $\mathcal D(w_j)$, which is a polynomial algebra in a differential operator $\delta_j$ (the Hermite, Laguerre, or Jacobi operator). 
Hence, our focus now shifts to understanding the structure of the modules \(D(w_j, w_k)\) when \(j \neq k\).
The main result  is that  there exists a differential operator $\mathcal T=\mathcal T_{j,k}$, with polynomial coefficients such that $$\mathcal{D}(w_{j}, w_{k})= \mathcal{T}\cdot \mathcal{D}(w_{k}).$$
 In other words, $\mathcal D(w_j,w_k)$ is a cyclic right module over the algebra $\mathcal{D}(w_{k})$, and we will provide the explicit expression of these operators $\mathcal T$ in all possible cases of Hermite, Laguerre, and Jacobi weights.

\begin{prop}\label{modules 0}
 If $w_{1}(x)w_{2}(x)^{-1}$ is not a rational function then $\mathcal{D}(w_{1},w_{2}) = \mathcal{D}(w_{2},w_{1})=0$.
\end{prop}
\begin{proof}
   From 
 Proposition \ref{diag D}, the algebra $\mathcal D(w_1\oplus w_2)$ is a diagonal algebra. Thus the module  $\mathcal D(w_1,w_2)$ and $ \mathcal{D}(w_{2},w_{1})$ are zero. 
\end{proof}

We have a nice condition for the module $D(w_1,w_2)$ to be zero in terms of the Darboux transformation of the weights.

\begin{thm}\label{Dw_1w_2zero}
    A classical scalar weight $w_{1}(x)$ is a Darboux transformation of a classical scalar weight $w_{2}(x)$ if and only if $\mathcal{D}(w_{1},w_{2}) \neq 0$.
\end{thm}

\begin{proof}
If $w_1$ is a Darboux transformation of a weight  $w_2$ then 
%from Proposition \ref{art} 
there exists a nonzero differential operators 
$\mathcal{V}$ and $\mathcal{N}$ such that 
$\left(\begin{smallmatrix} 0 && \mathcal{V} \\ \mathcal{N} && 0 \end{smallmatrix}\right) \in \mathcal{D}(w_1 \oplus w_2)$. Therefore $\mathcal{V}$ belongs to $D(w_1,w_2)$ and $\mathcal{N}$ belongs to $D(w_2,w_1)$.

\smallskip

We assume that $w_1$ is not a Darboux transformation of $w_2$ and we consider the cases of Hermite, Laguerre, and Jacobi weights separately.  
For different shifted Hermite weights, we always have that 
$w_1(x)w_2(x)^{-1}$ 
is not a rational function, and hence $\mathcal D(w_1,w_2)=0$.

If a Laguerre weight $w_1$ is not a Darboux transformation of $w_2$ 
%are Laguerre weights,  $w_1\not \sim w_2$, 
then, from Theorem \ref{laguerre scalar equivalences}, we have that the module satisfies $\mathcal D(w_1,w_2)=0$.

 Now, consider Jacobi weights $w_1= w_{\alpha, \beta}$ and $w_2= w_{\widetilde\alpha, \widetilde\beta}$,  such that $w_1\not \sim w_2$. Again, if $\alpha-
 \widetilde \alpha$ or $\beta-\widetilde\beta \not \in \mathbb Z$ then  $w_1(x)w_2(x)^{-1}$ is not a rational function, and hence, $\mathcal D(w_1,w_2)=0$. 
 From Theorem \ref{Jacobi} we only have to consider the case $\widetilde{\alpha} = \alpha + m+k$ and  $\widetilde{\beta} =\beta - k$, for some $k,m\in \mathbb N$.  
    In this case $w_{\widetilde\alpha,\widetilde\beta} = w_{{\alpha}+k+m,{\beta} -k}$ is a Darboux transformation of $ w_{{\alpha}+m,{\beta}}$. 
    In particular there exist degree-preserving differential operators $\mathcal{V} \in \mathcal{D}
    (w_{\widetilde{\alpha},\widetilde{\beta}}, w_{\alpha+m,\beta}) $ and $\mathcal{N} \in \mathcal{D}(w_{\alpha+m,\beta}, w_{\widetilde{\alpha},\widetilde{\beta}})$. 
    
If there exists a nonzero differential operator $0\neq \mathcal{S} \in \mathcal{D}(w_{\alpha,\beta}, w_{\widetilde{\alpha},\widetilde{\beta}})$ then, by  Proposition \ref{alg 22}, there also exists 
    $0\neq \mathcal{T} \in \mathcal{D}( w_{\widetilde{\alpha},\widetilde{\beta}}, w_{\alpha,\beta})$.
Thus  we get 
     $U_1= \mathcal{S}  \mathcal{V}\in \mathcal{D}(w_{{\alpha},{\beta}},w_{{\alpha}+m,{\beta}})$, $ 
     U_2= \mathcal{N}  \mathcal{T}\in \mathcal{D}(w_{{\alpha+m},{\beta}},w_{{\alpha},{\beta}})$. We then have   
     $$ 
     U_1U_2\in \mathcal{D}(w_{\alpha,\beta}), \quad 
  \text{ and }    \qquad U_2U_1\in \mathcal{D}(w_{\alpha+m,\beta})
     .$$ 
By proceeding in the same way as in Proposition \ref{pop}, we obtain a contradiction and therefore $\mathcal{D}(w_{\alpha,\beta},w_{\widetilde{\alpha},\widetilde{\beta}}) = 0$, which concludes the proof.  
\end{proof}

To study the structure of $\mathcal D(w_1, w_2)$ as a right $\mathcal D(w_2)$-module when $w_1$ is a Darboux transformation of $w_2$ we only have to consider the following situations:
\begin{align*}
\text{ Laguerre weights } \qquad\quad w_1&=w_\alpha, \quad w_2= w_{\alpha +k}, \quad \text { for } k\in \mathbb Z,\\
\text{  Jacobi weights } \qquad\quad w_1&=w_{\alpha,\beta +k}, \quad w_2= w_{\alpha+k,\beta }, \quad \text { for } k\in \mathbb Z.   
\end{align*}

We will see that the description of these modules is closely related to the Darboux transformation of weights. 
We will take full advantage of the Darboux transformations given in Section \ref{classical-section} to obtain an explicit generator of the right-module $\mathcal{D}(w_1,w_2)$ over the algebra $\mathcal D(w_2)$.

\smallskip

For  {\em Laguerre  weights}  
  $w_{\alpha}(x) = e^{-x}x^{\alpha}$,  
 we have proved  in Theorem \ref{laguerre scalar equivalences} that  
$\mathcal{V} = \partial - 1 \in \mathcal{D}(w_{\alpha},w_{\alpha+1})$ and $\mathcal{N} = \partial x + \alpha + 1 \in \mathcal{D}(w_{\alpha+1},w_{\alpha}).$ 
For $k \in \mathbb N$, we define the operator \begin{equation} \label{N-Lag-def} \mathcal{N}^{(k)} = \mathcal{N}_{k}\mathcal{N}_{k-1}\cdots \mathcal{N}_{1},  
\quad \text{ where }  
\quad \mathcal{N}_{j} = \partial x + \alpha + j.
\end{equation}
Hence, 
$$\mathcal{V}^{k}  = (\partial - 1)^{k} 
\in \mathcal{D}(w_{\alpha},w_{\alpha+k}) \quad \text{ and } \quad \mathcal{N}^{(k)} \in \mathcal{D}(w_{\alpha+k},w_{\alpha}).$$

The differential operator  $\mathcal{V}^{k}$ is of order $k$ with leading coefficient $1$, and $\mathcal{N}^{(k)}$ is a differential operator of order $k$ with leading coefficient $x^{k}$. We also have that 
\begin{equation}\label{lag pr}\mathcal{V}^{k}\cdot \mathcal{D}(w_{\alpha+k},w_{\alpha}) \subseteq \mathcal{D}(w_{\alpha}) \qquad \text{and} \qquad \mathcal{N}^{(k)}\cdot \mathcal{D}(w_{\alpha},w_{\alpha+k}) \subseteq \mathcal{D}(w_{\alpha+k}).
\end{equation}

\smallskip
\begin{prop} \label{l2}
    Let $w_{\alpha}(x) = e^{-x}x^{\alpha}$ be a Laguerre weight and $k \in \mathbb{N}$. 
    Then 
    $$\mathcal{D}(w_{\alpha},w_{\alpha+k}) = \mathcal{V}^{k} \cdot \mathcal{D}(w_{\alpha+k})\quad \text{ and } \quad \mathcal{D}(w_{\alpha+k},w_{\alpha}) = \mathcal{N}^{(k)}\cdot  \mathcal{D}(w_{\alpha}), $$
where $\mathcal{V} = \partial - 1$ and $\mathcal{N}^{(k)}$ are defined as in \eqref{N-Lag-def}.
\end{prop}
\begin{proof}
     Let us recall that $\mathcal{D}(w_{\alpha}) = \mathbb{C}[\delta_{\alpha}]$, where $\delta_{\alpha} = \partial^{2}x + \partial (\alpha+1-x)$. % In particular the leading coefficient of $(\delta_{\alpha})^{n}$ is $x^{n}$. 
To see that $\mathcal{D}(w_{\alpha},w_{\alpha+k}) \subseteq \mathcal{V}^{k}\cdot \mathcal{D}(w_{\alpha+k})$, we take $\mathcal{A} \in \mathcal{D}(w_{\alpha},w_{\alpha+k})$, $\mathcal A\neq 0$.  Since $\mathcal{N}^{(k)} \in \mathcal{D}(w_{\alpha+k},w_{\alpha})$, we have that 
    $\mathcal{A}\mathcal{N}^{(k)}\in \mathcal{D}(w_{\alpha})$. 
     Therefore $\mathcal{A}\mathcal{N}^{(k)} = p(\delta_{\alpha})$ for some polynomial $p\in \mathbb{C}[x]$ of degree $m$. 
     %The leading coefficient of $\mathcal{A}\mathcal{N}^{(k)}$ is $cx^{m}$ for some $c\in \mathbb{C}$ and 
     Then 
     the leading coefficient of $\mathcal{A}$ is $cx^{m-k}$, for some $c\in \mathbb{C}$  %in particular  $m \geq k$. 
     %The order of the differential operator $p(\delta_{\alpha})$ is $2m$, therefore 
     and the order of % the nonzero operator 
     $\mathcal{A}$ is $\operatorname{ord}(\mathcal A)=2m - k \geq k$. 
     
     We proceed by induction on the order of $\mathcal{A}$. For $m = k$, it follows that the order of $\mathcal{A}$ is $k$ and its leading coefficient is $c$ for some $c\in \mathbb{C}$. Thus $\mathcal{A} - c\mathcal{V}^{k}\in \mathcal{D}(w_{\alpha},w_{\alpha + k})$ and $\operatorname{ord}(A- c\mathcal{V}^{k})<k$. Hence $\mathcal{A} = c\mathcal{V}^{k}$ and the statement holds. 
     Now, suppose that $\operatorname{ord}(A)>m \geq k$. Therefore $\operatorname{ord}\big (\mathcal{A} - c\mathcal{V}^{k}(\delta_{\alpha+k})^{m-k} \big ) < 2m-k$. 
     Thus by  inductive hypothesis $\mathcal{A} -  c \mathcal{V}^{k}(\delta_{\alpha+k})^{m-k} \in  \mathcal{V}^{k}\mathcal D(w_{\alpha+k})$, 
     %for some $D\in \mathcal{D}(w_{\alpha+k})$. Hence $\mathcal{A} = \mathcal{V}^{k}(D+(\delta_{\alpha+k})^{m-k}) \in \mathcal{V}^{k}\cdot \mathcal{D}(w_{\alpha+k}) $,
     which concludes the proof.  

     The second part of the statement follows similarly.
\end{proof}

\ 

Now, we consider  Jacobi weights
    $w_{\alpha,\beta}(x) = (1-x)^{\alpha}(1+x)^{\beta}$, for  $\alpha,\, \beta > -1$ and $x\in (-1.1)$.
    We have proved  in Theorem \ref{Jacobi} that 
    $$\mathcal{V} = \partial (1+x) + \beta + 1 \in \mathcal{D}(w_{\alpha,\beta+1},w_{\alpha+1,\beta}) \quad \text{ and } \quad \mathcal{N} = \partial (x-1) + \alpha + 1 \in \mathcal{D}(w_{\alpha+1,\beta},w_{\alpha,\beta+1}). $$

    For $k\in \mathbb N$, we define the operators 
    \begin{equation}\label{NV-op-Jaco}
    \begin{split}
      \mathcal{V}^{(k)} &= \mathcal{V}_{k} \, \mathcal{V}_{k-1}\cdots \mathcal{V}_{1} \quad \text{ where }\quad \mathcal{V}_{j} = \partial(1+x) + \beta + j, 
      \\
      \mathcal{N}^{(k)} &= \mathcal{N}_{k} \, \mathcal{N}_{k-1} \cdots \mathcal{N}_{1}   \quad \text{ where } \quad \mathcal{N}_{j} = \partial(x-1) + \alpha + j.
    \end{split}
    \end{equation}

    As a consequence we have that $\mathcal{V}^{(k)}  \in \mathcal{D}(w_{\alpha,\beta+k},w_{\alpha+k,\beta})$ and $\mathcal{N}^{(k)}  \in \mathcal{D}(w_{\alpha+k,\beta},w_{\alpha,\beta+k})$. 
    The operators $\mathcal{V}^{(k)}$ and $\mathcal{N}^{(k)}$ are differential operators of order $k$ with leading coefficients $(1+x)^{k}$ and $(x-1)^{k}$, respectively. We also have that $$\mathcal{V}^{(k)} \cdot \mathcal{D}(w_{\alpha+k,\beta},w_{\alpha,\beta+k}) \subseteq \mathcal{D}(w_{\alpha,\beta+k}) \quad \text { and } \quad  \mathcal{N}^{(k)} \cdot \mathcal D(w_{\alpha,\beta+k},w_{\alpha+k,\beta}) \subseteq \mathcal{D}(w_{\alpha+k,\beta}).$$

\smallskip
\begin{prop} \label{jac 1} For Jacobi weights $w_{\alpha, \beta}$, $\alpha,\beta > -1$, and $k\in \mathbb N$ we have
    $$\mathcal{D}(w_{\alpha,\beta+k},w_{\alpha+k,\beta}) = \mathcal{V}^{(k)} \cdot \mathcal{D}(w_{\alpha+k,\beta}) \quad  \text{ and } \quad \mathcal{D}(w_{\alpha+k,\beta},w_{\alpha,\beta+k}) = \mathcal{N}^{(k)} \cdot \mathcal{D}(w_{\alpha,\beta+k}), $$
    where $\mathcal V $ and $\mathcal{N }$ are defined in \eqref{NV-op-Jaco}.
\end{prop}
\begin{proof}
    The proof of this proposition is analogous to the proof of Proposition \ref{l2}.
\end{proof}

\ 

We can summarize all these results in the following theorem.

\begin{thm}\label{generatorT}
    Given  $w_1, w_2$ two classical scalar weights, the right-module $\mathcal D(w_1,w_2)$ is cyclic over the algebra $\mathcal{D}(w_{2})$, i.e., there exists a  differential operator $\mathcal T= \mathcal T_{w_1,w_2} \in \Omega[x] $ such that  $$\mathcal{D}(w_{1}, w_{2})= \mathcal{T}\cdot \mathcal{D}(w_{2}).$$
Moreover, 
\begin{enumerate}
    \item [(i) ]If $w_1$ is not a Darboux transformation of $w_2$, then $\mathcal T=0$.
\item [(ii) ]If $w_1=w_2$, then $\mathcal T=1$.
\item [(iii) ] For Laguerre weights $w_\alpha$,  $k\in \mathbb N$,  we have 
$$ \mathcal T_{w_\alpha, w_{\alpha+k}} = (\partial -1)^k,  \qquad  \mathcal T_{w_\alpha+k, w_{\alpha}}= \mathcal N^{(k)}, \quad \text{ with } \mathcal N^{(k) } \text{ as defined in } \eqref{N-Lag-def}.$$
\item [(iv) ] For Jacobi weights $w_{\alpha,\beta}$ we have 
$$ \mathcal T_{w_{\alpha,\beta+k}, w_{\alpha+k,\beta}} = \mathcal V^{(k)},  \qquad  \mathcal T_{w_{\alpha+k,\beta}, w_{\alpha,\beta+k}}= \mathcal N^{(k)}, \quad \text{ with }  \mathcal V^{(k)} , \mathcal N^{(k) } \text{ as defined in } \eqref{NV-op-Jaco}.$$
\end{enumerate}
\end{thm}

\medskip

In this way, we can determine the algebra $\mathcal D(W)$ for $W= w_{1}\oplus w_2 \oplus \cdots \oplus w_{N}$, a direct sum of classical scalar weights of the same type (shifted Hermite, Laguerre, or Jacobi). Recall that the algebra $\mathcal D(w_j)$ is a polynomial algebra in the differential operator $\delta_j=\delta$ given in Table \ref{weights-table}. 

\begin{thm}\label{algdirect}
    Let $W= w_{1}\oplus w_2\oplus  \cdots \oplus w_{N}$ be a direct sum of classical scalar weights of the same type (shifted Hermite, Laguerre, or Jacobi). Then  
    \begin{equation*}
    \mathcal{D}(w_{1}\oplus \cdots \oplus w_{N} ) = \sum_{i=1}^{N}\sum_{j=1}^{N} \mathcal{T}_{w_i,w_j}\cdot \mathcal{D}(w_{j})\, E_{i,j},
\end{equation*}
where the differential operators $\mathcal T_{w_i,w_j}  $ are those given in Theorem \ref{generatorT}.
\end{thm}

\begin{cor}\label{generators}
  Let $W= w_{1}\oplus w_2\oplus  \cdots \oplus w_{N}$ be a direct sum of classical scalar weights of the same type (shifted Hermite, Laguerre, or Jacobi). Then  the algebra $\mathcal D(W)$ is generated by the set 
  $$\left  \{  \delta_{w_j} E_{j,j}\, , \mathcal T_{w_i,w_j} E_{i,j}, E_{j,j} \, : \, 1\leq i\neq j\leq N\right \}.$$
\end{cor}

\medskip
\noindent {\bf  Example.} To illustrate these results, let us consider a direct sum of scalar Laguerre weights $ W=w_{\alpha-1 }\oplus w_{\alpha}$. 
We have that 
$$\mathcal{D}(w_{\alpha-1}\oplus w_{\alpha}) = \begin{pmatrix} \mathcal{D}(w_{\alpha-1}) && \mathcal{V} \cdot \mathcal{D}(w_{\alpha}) \\\mathcal{N} \cdot \mathcal{D}(w_{\alpha-1}) && \mathcal{D}(w_{\alpha}) \end{pmatrix} = \begin{pmatrix} 1 && \mathcal{V} \\ \mathcal{N} && 1 \end{pmatrix} \begin{pmatrix} \mathcal{D}(w_{\alpha-1}) && 0 \\0 && \mathcal{D}(w_{\alpha}) \end{pmatrix},$$
where $\mathcal{V} = \partial - 1$ and $\mathcal{N} = \partial x + \alpha $. Therefore the algebra $\mathcal D(w_{\alpha-1}\oplus w_{\alpha}) $ can be generated by  $\quad 
\left(\begin{smallmatrix}
    1& 0\\  0 & 0
\end{smallmatrix} \right)$, $\left(\begin{smallmatrix}
    0& 0\\  0 & 1
\end{smallmatrix} \right)$   and 

$$  \left(\begin{matrix}
    \partial^2 x+ \partial (\alpha -x) & 0\\0 &0
\end{matrix} \right), \quad \left(\begin{matrix}
    0& 0\\0 & \partial^2 x+ \partial (\alpha+1 -x) \end{matrix} \right) , \quad 
\left(\begin{matrix}
    0& \partial -1\\0 & 0
\end{matrix} \right), \quad \left(\begin{matrix}
    0& 0\\  \partial x + \alpha  & 0
\end{matrix} \right). 
$$ 

\ 

Now, we obtain an important consequence of Theorem \ref{algdirect}.  In particular, we use this result to complete the proof of Theorem \ref{com 2}.

\begin{prop}\label{com}
 Let $W= w_{1}\oplus \cdots \oplus w_{N}$ be a direct sum of classical scalar weights, then the algebra $\mathcal D (W)$ is commutative if and only if $w_{i}$ is not a Darboux transformation of $ w_{j}$ for all $i\not= j$.
\end{prop}
\begin{proof}
    If $w_{i} \sim w_{j}$ for some $i\not= j$, then the nonzero differential operator $\mathcal{T}E_{i,j}$ belongs to $\mathcal{D}(W)$ and it does not commute with $E_{ii}$. Thus, $\mathcal{D}(W)$ is not a commutative algebra. Now, assume that $w_{i}\not \sim w_{j}$ for all $i\not = j$, then $\mathcal{D}(W) = \mathcal{D}(w_{1}) \oplus \cdots \oplus \mathcal{D}(w_{N})$ which is a commutative algebra. 
\end{proof}

\smallskip

We can use the knowledge of the algebra $\mathcal{D}(W)$ for a diagonal weight $W$ to determine when two diagonal weights are Darboux-equivalent. As an initial step, we establish a necessary condition for this equivalence.

\begin{prop}\label{mismos pesos}
   Let $W= w_{1}\oplus  \cdots  \oplus w_{N}$ and $\widetilde W= \widetilde w_{1} \oplus  \cdots \oplus \widetilde w_{N}$ be direct sums of scalar weights. 
    If $W$ is a Darboux transformation of $\widetilde W$, then each weight $ w_i$ is a Darboux transformation of some weight $\widetilde w_j$. 
\end{prop}
\begin{proof} 
Without loss of generality, we can assume that $w_1$ is not a Darboux transformation of any weight $\widetilde w_j$, for all $j = 1,\ldots, N$. Thus $\mathcal D(w_1, \widetilde w_j)=0$. 
   Let $P_n=\operatorname{diag}(p_n^{w_1}, \dots ,  p_n^{w_1} )$ and $\widetilde P_n=\operatorname{diag}(\widetilde p_n^{\widetilde w_1}, \dots ,  \widetilde p_n^{\widetilde w_1} )$  be the monic sequence of orthogonal polynomials associated to $W$ and $\widetilde W$. 
   If $W$ is a Darboux transformation of $\widetilde W$, then
   there exists a degree preserving differential operator $\mathcal V$ such that $P_n\cdot \mathcal V= A_n \widetilde P_n$, for some matrix $A_n$. Hence 
   $$ p_n\cdot \mathcal V_{1,j} = (A_n)_{1,j}\,  \widetilde p_n^{\widetilde w_j}, \qquad \text { for } j=1,\cdots N. $$  
   This implies that 
   the operator $\mathcal V_{1,j}\in \mathcal{D}(w_1, \widetilde w_j)=0$, i.e., the first row of $\mathcal V$ is zero. This is a contradiction because  $\mathcal V$ is a degree-preserving operator.
\end{proof}

\begin{remark} 
The condition stated in Proposition \ref{mismos pesos} is not sufficient to ensure that  $W = w_{1} \oplus \cdots \oplus w_{N}$ is Darboux-equivalent  to  $\widetilde{W} = \widetilde{w}_{1} \oplus \cdots \oplus \widetilde{w}_{N}$. For example, if $w_1\not\sim w_2$ then
$$W=\left(\begin{smallmatrix}
    w_1& & &\\ & w_1 & & \\ &&w_1 & \\ &&& w_2  
\end{smallmatrix} \right) \quad \text{ is not a Darboux transformation of } \quad 
\widetilde W=\left(\begin{smallmatrix}
    w_1& & &\\ & w_1 & & \\ &&w_2 & \\ &&& w_2  
\end{smallmatrix} \right) . $$
In fact, if $W$ were a Darboux transformation of $\widetilde W$ then there would exist a degree preserving differential operator $\mathcal V$ such that $P_n\cdot \mathcal V= A_n \widetilde P_n$, for some nonsingular matrix $A_n$.
By proceeding in the same way 
as in the proof of  Proposition \ref{mismos pesos}, we obtain that $\mathcal V$  and the matrices $A_n$ are of the form
$$\mathcal V= \left(\begin{smallmatrix}
    \mathcal V_{11} & \mathcal V_{12} & 0 & 0 \\ \mathcal V_{21} &  \mathcal V_{22} &0 &0 \\ \mathcal V_{31} &  \mathcal V_{32} &0 &0 \\ 0 & 0&  \mathcal V_{43} & \mathcal V_{44}  \end{smallmatrix} \right),  \qquad  \quad 
    A_n= \left(\begin{smallmatrix}
   a_{11} & a_{12} & 0 & 0 \\ a_{21} &  a_{22} &0 &0 \\ a_{31} &  a_{32} &0 &0 \\ 0 & 0&  a_{43} & a_{44}  \end{smallmatrix} \right), $$ 
   with 
   $ \mathcal V_{j1}, \mathcal V_{j2} \in \mathcal D(w_1) $,  ($j=1,2,3$)   and $ \mathcal V_{43}, \mathcal V_{44} \in \mathcal D(w_2)$. 
But these matrices $A_n$ are all singular matrices.
\end{remark}

Finally, we can determine when two diagonal weights are obtained  from each ather,  through a Darboux transformation.

\begin{thm} \label{Darboux-direct sum}
   Let $W= w_{1}\oplus \cdots \oplus w_{N}$ and $\widetilde W= \widetilde w_{1}\oplus \cdots \oplus \widetilde w_{N}$ be direct sums of classical scalar weights of the same type. 
    Then $\widetilde W$ is a Darboux transformation of $W$ if and only if there exists a permutation $\sigma$ such that $\widetilde w_{j} $ is a Darboux transformation of ${w}_{\sigma(j)}$, for all $j = 1,\ldots, N$.
\end{thm}

\begin{proof}
If the set of weights $\{ \widetilde w_1, \dots, \widetilde w_N\}$ is a permutation of  $\{ w_1, \dots, w_N\}$, then it is clear that $W$ is a Darboux transformation of $\widetilde W$. 
On the other hand, from Proposition \ref{mismos pesos}, we can assume that 
   \begin{align*}
       W & = w_{1}I_{r_{1}\times r_{1}
    } \oplus w_{2}I_{r_{2}\times r_{2}} \oplus \cdots \oplus w_{m} I_{r_{m}\times r_{m}},   
 \qquad 
    \widetilde{W}  = w_{1}I_{s_{1}\times s_{1}} \oplus w_{2}I_{s_{2}\times s_{2}} \oplus \cdots \oplus w_{m} I_{s_{m}\times s_{m}}, 
    \end{align*}
 and the weights $w_{i}$ are not a Darboux transformation of $ w_{j}$ for every $j \neq  i$.    
Note that in particular, $r_{1} + \cdots + r_{m} = N = s_{1} + \cdots + s_{m}$ and  $\mathcal D(w_i,w_j)=0$, for $i\neq j$ (see Proposition \ref{generatorT}).

Since $\widetilde W$ is a Darboux transformation of $W$ 
then there exists a degree preserving operator $\mathcal V$ such that $P_n \cdot \mathcal V= A_n \widetilde P_n$, where 
$$P_n(x)=\left(\begin{smallmatrix}p_{n}^{w_{1}} I_{r_{1} \times r_{1}} && && \\
    && \ddots && \\
    && && p_{n}^{w_{m}}I_{r_{m}\times r_{m}}\end{smallmatrix} \right )\; \text{ and } \; 
    \widetilde P_n(x)= \left(\begin{smallmatrix} p_{n}^{w_{1}} I_{s_{1} \times s_{1}} && && \\
    && \ddots && \\
    && && p_{n}^{w_{m}}I_{s_{m}\times s_{m}} \end{smallmatrix} \right ) $$
are the monic sequence of orthogonal polynomials associated with $W$ and $\widetilde W$, respectively.
Then we have that $\mathcal V$ is a block diagonal matrix
   $$\mathcal V= \mathcal V_1\oplus \cdots \oplus \mathcal V_m \qquad  \text{ with }  \quad \mathcal V_j\in \operatorname{Mat}_{r_j\times s_j} (\mathcal D(w_j)). $$
    Therefore, the matrices $A_n$ are also block matrices 
    $ A_n= A_1(n) \oplus \cdots \oplus A_m(n)$ with   $A_j(n)\in \operatorname{Mat}_{r_j\times s_j}(\mathbb C).$

Now, if $r_{i} \not= s_{i}$ for some $i$, it follows that the columns or the rows of $A_{n}$ must be linearly dependent, which means that $A_{n}$ is a singular matrix and this is a contradiction because $\mathcal{V}$ is a degree-preserving differential operator. 
 Therefore  $r_{i} = s_{i}$ for all $i$ and the proof is now complete. 
\end{proof}

\section{An example arising from spherical functions}
\label{Sect-sph}
The harmonic analysis on compact symmetric spaces has played a fundamental role in the construction of families of matrix-valued orthogonal polynomials which are eigenfunctions of a second-order differential operator. 
A first example was given in \cite{GPT02} for the symmetric pair $(\mathrm{SU}(3), \mathrm{U}(2) )$.   
The following example, given in \cite{PT07} for $N=2$,  arises by studying matrix-valued spherical functions associated with the pair $(G,K)$ with $G=\mathrm{SU}(m+1)$ and $ K=\mathrm{U}(m)$. 

After the change of variable $x=2t-1$, the weight matrix $W$ is given by 
\begin{equation}\label{Jacobi-2x2}
    W(x) = (1-x)^{\alpha}(1+x)^{\beta}\begin{pmatrix} 4(\alpha+1) - 2k( 1-x) && 2(\alpha + 1 - k)(1-x) \\ 2(\alpha +1 -k)(1-x) && (\alpha +1 -k)(1-x)^{2}\end{pmatrix},
\end{equation}
where  $\alpha, \beta> -1$, $0<k<\beta+1$ and $x\in (-1,1)$. 
This is an irreducible weight and the algebra $\mathcal{D}(W)$ contains the second-order differential operator
\begin{align*} 
D_{1} & = \partial^{2} (1-x^{2})I + \partial \begin{pmatrix}\beta - \alpha +1 -x(\alpha + \beta +3) && -2 \\ 0 && \beta - \alpha - 2 - x(\alpha+\beta+4)  \end{pmatrix}\\
        & \quad + \begin{pmatrix}\alpha+\beta-k+2 && 0 \\ \alpha + 1 - k && 0 \end{pmatrix}. 
        \end{align*} 
The $W$-symmetric differential operator       
\begin{equation*}
    \begin{split}
        D_{2} & = \partial^{2} \begin{pmatrix} 1 - x^{2} && -2(x+1) \\ 0 && 0 \end{pmatrix} + \partial \begin{pmatrix} \alpha + \beta -2k +3 -x(\alpha + \beta + 3) && -2(\alpha+\beta - k +3) \\ (1-x)(\alpha+1-k) && -2(\alpha +1 -k) \end{pmatrix} \\
        & \qquad + \begin{pmatrix}-k(\alpha + \beta + 2- k) && 0 \\ -k(\alpha + 1 -k) && 0 \end{pmatrix}
    \end{split}
\end{equation*}
also belongs to $\mathcal D(W)$ and commutes with $D_1$. 

In \cite{PR08}, a sequence of matrix-valued orthogonal polynomials is given in terms of the matrix-valued hypergeometric functions. See also \cite{P08} for  more explicit expressions of these polynomials 
in the case $\beta=1$.

\smallskip
Now we prove that this irreducible weight is a Darboux transformation of $w_{\alpha+1,\beta}\oplus w_{\alpha+1, \beta+1}$, a direct sum of classical Jacobi weights. 

\begin{prop} \label{Jacobi2x2-Darboux} 
The  weight matrix
$W(x)$ is a Darboux transformation of the direct sum 
$$w(x) = \begin{pmatrix} (1-x)^{\alpha+1}(1+x)^{\beta} && 0 \\ 0 && (1-x)^{\alpha+1}(1+x)^{\beta+1}\end{pmatrix}.$$
\end{prop}

\begin{proof}
We take the differential operator  $D \in \mathcal{D}(w)$,
\begin{equation*}
    \begin{split}
        D & = \partial^{2} \begin{pmatrix} 2k(x^{2}-1) && 0 \\ 0 && (\alpha+1-k)(x^{2}-1)\end{pmatrix} \\
        & \quad + \partial \begin{pmatrix} 2k(\alpha-\beta+1+x(\alpha+\beta+3)) && 0 \\ 0 && (\alpha+1-k)(\alpha-\beta+x(\alpha+\beta+4))\end{pmatrix} \\ 
        & \quad + \begin{pmatrix} 2k^{2}(\alpha + \beta - k + 2) && 0 \\ 0 && (k+1)(\alpha+1-k)(\alpha+\beta-k+2) \end{pmatrix}.
 \end{split}
\end{equation*}

 We factorize $D$ as $D = \mathcal{V}\mathcal{N}$, where
 
$$\mathcal{V} = \partial \begin{pmatrix}x -1 && 2 \\ 0 && 1+x \end{pmatrix} + \begin{pmatrix}k && 0 \\ k - \alpha - 1 && \alpha + \beta - k + 2 \end{pmatrix} $$
and 
$$\mathcal{N} = \partial \begin{pmatrix}2k(x+1) && 2(k-\alpha-1) \\ 0 && (\alpha + 1 -k)(x-1) \end{pmatrix} + \begin{pmatrix}2k(\alpha + \beta - k + 2) && 0 \\ 2k(\alpha+1-k) && (k+1)(\alpha+1-k) \end{pmatrix}. $$

\smallskip
Let $J_{n}^{(\alpha,\beta)}(x)$ be the sequence of monic orthogonal polynomials for the scalar Jacobi weight $w_{\alpha,\beta}(x) = (1-x)^{\alpha}(1+x)^{\beta}$. It follows that $Q_{n}(x) = \begin{pmatrix} J_{n}^{(\alpha+1,\beta)}(x) && 0 \\ 0 && J_{n}^{(\alpha+1,\beta+1)}(x) \end{pmatrix}\cdot \mathcal{V}$ is a sequence of orthogonal polynomials for $W$. Hence $W$ is a Darboux transformation of $w$. 
\end{proof}

\begin{remark}
The differential operator $\widetilde{D} = \mathcal{N}\mathcal{V}$ belongs to $\mathcal{D}(W)$ and it is given by  
$$\widetilde{D} = -(\alpha + 1 -k)D_{1} +(\alpha-3k+1)D_{2} + (k+1)(\alpha+1-k)(\alpha+\beta-k+2)I.$$
   
\end{remark}

\medskip
Now, by using the techniques and concepts developed in this paper, we can obtain new results about the structure of algebra $\mathcal D(W)$.

\begin{prop} Let $W$ be the   weight matrix  
introduced in \eqref{Jacobi-2x2}. Then the following holds:
\begin{enumerate}
    \item [i)]  The algebra 
    $\mathcal{D}(W)$ is a commutative algebra.
    \item [ii)]  
    $\mathcal{D}(W)$ does not contain odd-order differential operators and it is generated by $\{I, D_{1}, D_{2} \}$.

\end{enumerate}

\end{prop}
\begin{proof}  \emph{ i)}  
The weight matrix  $W(x)$ is a Darboux transformation of the diagonal weight $w(x)=w_{\alpha+1,\beta}\oplus w_{\alpha+1,\beta+1}$.  
From Theorem \ref{Jacobi},
the scalar weight $w_{\alpha+1,\beta}$ is not a Darboux transformation of $w_{\alpha+1,\beta+1}$, hence $\mathcal D(w)$ is a commutative algebra by   Proposition  \ref{com}.
Now, the statement follows from Theorem \ref{com 2}. 

    \emph{ ii)} By Theorem \ref{algdirect}, it follows that $\mathcal{D}(w) = \mathcal{D}(w_{\alpha+1,\beta}) \oplus \mathcal{D}(w_{\alpha+1,\beta+1})$. If $\mathcal{A} \in \mathcal{D}(W)$, with  $\operatorname{Ord}(\mathcal{A}) = 2m+1$, then the differential operator $\mathcal{B}=\mathcal{V}\mathcal{A}\mathcal{N}$ belongs to $\mathcal{D}(w)$, by Proposition \ref{alg debil}.   
    Since the leading coefficients of $\mathcal{V}$ and $\mathcal{N}$ are not singular, it follows that $\operatorname{Ord}(\mathcal{B}) = 2m+3$, but $\mathcal{D}(w)$ has only even-order differential operators which leads to a contradiction.
    
    To prove that every differential operator $\mathcal{A} \in \mathcal{D}(W)$ is generated by $\{I,D_{1},D_{2} \}$, we will proceed by induction on the order of $\mathcal{A}$. First, by the equations of symmetry, it is easy to check that the statement holds for differential operators of order less than or equal to $2$. Now, let $m \geq 2$, given a differential operator $\mathcal{A} = \sum_{j=0}^{2m} \partial^{j}F_{j} \in \mathcal{D}(W)$, we have that the differential operator $\mathcal{B} = \mathcal{V}\mathcal{A}\mathcal{N}$ of order $2m+2$ belongs to $\mathcal{D}(w)$. Therefore its leading coefficient is of the form $\begin{pmatrix} \ell_{1} (1-x^{2})^{m+1} && 0 \\ 0 && \ell_{2}(1-x^{2})^{m+1} \end{pmatrix}$ for some complex numbers $\ell_{1}$, $\ell_{2}$. On the other hand, the leading coefficient of $\mathcal{B}$ is $\begin{pmatrix}x -1 && 2 \\ 0 && 1+x \end{pmatrix}F_{2m}\begin{pmatrix}2k(x+1) && 2(k-\alpha-1) \\ 0 && (\alpha + 1 -k)(x-1) \end{pmatrix}$. From here, we obtain that the leading coefficient of $\mathcal{A}$ is 
    $$F_{2m}(x) = \begin{pmatrix}-\frac{\ell_{1}}{2k} (1-x^{2})^{m} && 2(-\frac{\ell_{1}}{2k}+\frac{\ell_{2}}{(\alpha+1-k)})(1+x)(1-x^{2})^{m-1} \\ 0 && -\frac{\ell_{2}}{(\alpha+1-k)}(1-x^{2})^{m} \end{pmatrix}.$$
    Thus, the differential operator $\mathcal{A} - \left(-\frac{\ell_{2}}{(\alpha+1-k)}D_{1}^{m} + \left(-\frac{\ell_{1}}{2k} + \frac{\ell_{2}}{(\alpha+1-k)} \right) D_{2}D_{1}^{m-1} \right)$ belongs to $\mathcal{D}(W)$ and it is of order less than $2m$. Hence, by our inductive assumptions, it follows that the algebra $\mathcal{D}(W)$ is generated by $\{I, D_{1}, D_{2}\}$.
\end{proof}

\ 

Related to the  solutions of the Matrix Bochner Problem given in \cite{CY18}, we recall that given a weight matrix $W$ of size $N\times N$, we said that the algebra $\mathcal{D}(W)$ is full if there exist nonzero $W$-symmetric operators $\mathfrak{D}_{1},\ldots,\mathfrak{D}_{N}$ in $\mathcal{D}(W)$,  such that  $$\mathfrak{D}_{i}\mathfrak{D}_{j} = 0 \text{ for $i\not= j\quad $ with } \quad \mathfrak{D}_{1} + \cdots + \mathfrak{D}_{N} \in \mathcal{Z}(W)$$ 
which is not a zero divisor. 

In the case of our weight $W$, given in \eqref{Jacobi-2x2},  we have the following result. 

\begin{prop}\label{fullfull}
    The algebra $\mathcal{D}(W)$ is a commutative full algebra.
\end{prop}
\begin{proof}
    The differential operators $\mathfrak{D}_{1} = D_{1}-D_{2}-(\alpha+\beta+2-k)(1+k)I$ and $\mathfrak{D}_{2} = D_{2}$ satisfy that 
    $$\mathfrak{D}_{1}\mathfrak{D}_{2} = 0 = \mathfrak{D}_{2}\mathfrak{D}_{1}$$ 
    and $\mathfrak{D}_{1} + \mathfrak{D}_{2}$ belongs to the center of $\mathcal{D}(W)$.
\end{proof}

\section{Some illustrative examples} \label{7}
In this section, we illustrate how the results obtained throughout this paper can be applied to specific examples. We will examine three concrete cases to demonstrate how these results provide new and detailed information about the irreducible weights $W$ and their algebra $\mathcal{D}(W)$.

\subsection{Example 1}

Let $w_{\alpha}(x) = e^{-x}x^{\alpha}$ be the classical scalar Laguerre weight with $\alpha > -1$. We  consider the diagonal weight 
$$W(x) = w_{\alpha+2}(x) \oplus w_{\alpha+1}(x) \oplus w_{\alpha}(x).$$
The sequence of monic orthogonal polynomials for the weight $W$ is given by 
$$P_{n}(x) = \left(\begin{smallmatrix} \ell_{n}^{(\alpha+2)}(x) & 0 & 0 \\ 0 & \ell_{n}^{(\alpha+1)}(x) & 0 \\ 0 & 0 & \ell_{n}^{(\alpha)}(x) \end{smallmatrix}\right),$$ where $\ell_{n}^{(\alpha)}(x)$ denotes the $n$-th monic orthogonal Laguerre polynomials  of parameter $\alpha$.

The algebra $\mathcal D(w_\alpha)$ associated to  the  scalar weight $w_\alpha$ is a polynomial algebra in the Laguerre differential operator $\delta_{\alpha} = \partial^{2}x + \partial (\alpha+1-x)$, i.e., $ \mathcal D(w_\alpha)= \mathbb{C}[\delta_{\alpha}]$, hence from Theorem \ref{algdirect} we known that the algebra $\mathcal{D}(W)$ is 
$$\mathcal{D}(W) = \begin{pmatrix} \mathbb{C}[\delta_{\alpha+2}] & (\partial x + \alpha + 2) \mathbb{C}[\delta_{\alpha+1}] & (\partial x + \alpha + 2)(\partial x + \alpha + 1) \mathbb{C}[\delta_{\alpha}] \\ (\partial - 1) \mathbb{C}[\delta_{\alpha+2}] & \mathbb{C}[\delta_{\alpha+1}] & (\partial x + \alpha + 1)\mathbb{C}[\delta_{\alpha}] \\ (\partial -1)^{2} \mathbb{C}[\delta_{\alpha+2}] & (\partial - 1 ) \mathbb{C}[\delta_{\alpha+1}] & \mathbb{C}[\delta_{\alpha}] \end{pmatrix}.$$

\smallskip
For  $a_{1}, \, a_{2} \in \mathbb{R}$, we consider  the sixth-order differential operator in the algebra $\mathcal{D}(W)$
$$D = \begin{psmallmatrix} -a_{1}^{2} \delta_{\alpha+2} + a_{1}^{2} + 1 & 0 & a_{2}a_{1}(\partial x + \alpha +2)(\partial x + \alpha + 1) (1-\delta_{\alpha}) \\ 0 & -a_{2}^{2}\delta_{\alpha+1}(\alpha+1-\delta_{\alpha+1})(\alpha-\delta_{\alpha+1})-a_{1}^{2}\delta_{\alpha+1}+1 & 0 \\ a_{2}a_{1}(\partial-1)^{2}(1-\delta_{\alpha+2}) & 0 & a_{2}^{2}(1-\delta_{\alpha})(\alpha+2-\delta_{\alpha})(\alpha+1-\delta_{\alpha})\end{psmallmatrix}. $$
We can verify that $D$  can be factorized as $D = \mathcal{V}\mathcal{N}$, with $\mathcal{V}$ and $\mathcal{N}$ the following degree-preserving differential operators 
\begin{align*}
        \mathcal{V} & = \partial^{3}\begin{psmallmatrix} 0 & 0 & 0 \\ 0 & a_{2}^{2}x^{3} & -a_{2}x^{2} \\ 0 & -a_{2} x & 0\end{psmallmatrix} + \partial^{2} \begin{psmallmatrix} 0 & 0 & 0 \\ 0 & 2a_{2}^{2}x^{2}(\alpha+2) & -2a_{2}x(\alpha+2) \\ 0 & -a_{2}(\alpha-3x+2) & 0 \end{psmallmatrix} + \partial \begin{psmallmatrix} 0 & -a_{1}x & 0 \\ -a_{1} & x(a_{2}^{2}(\alpha+2)(\alpha+1)+a_{1}^{2}) & -a_{2}(\alpha+2)(\alpha+1) \\ 0 & a_{2}(2 \alpha -3x +4) & 0 \end{psmallmatrix} \\
        & \quad + \begin{psmallmatrix}1 & -a_{1}(\alpha+2) & 0 \\ 0 & 1 & 0 \\ 0 & -a_{2}(\alpha+2) & 1 \end{psmallmatrix}, \\
     \intertext{ and }   \mathcal{N} & = \partial^{3} \begin{psmallmatrix}  0 & 0 & a_{1}a_{2}x^{3} \\ 0 & 0 & a_{2}x^{2} \\ 0 & a_{2}x & a_{2}^{2}x^{3}  \end{psmallmatrix} + \partial^{2} \begin{psmallmatrix} 0 & 0 & a_{1}a_{2}x^{2}(2\alpha + 7) \\ 0 & 0 & 2a_{2}x(\alpha+2) \\ 0 & a_{2}(\alpha - 3x +2) & a_{2}{2} x^{2} (2\alpha + 7) \end{psmallmatrix} + \partial \begin{psmallmatrix} a_{1}^{2}x & a_{1}x & a_{1}a_{2}(\alpha+5)(\alpha+2) \\ a_{1} & 0 & a_{2}(\alpha+2)(\alpha+1) \\ a_{1}a_{2}x & -a_{2}(2\alpha-3x+4) & a_{2}^{2}x(\alpha+5)(\alpha+2) \end{psmallmatrix} \\
        & \quad + \begin{psmallmatrix} a_{1}^{2}+1 & a_{1}(\alpha+2) & a_{1}a_{2}(\alpha+2)(\alpha+1) \\ 0 & 1 & 0 \\ a_{1}a_{2} & a_{2}(\alpha+2) & a_{2}^{2}(\alpha+2)(\alpha+1) \end{psmallmatrix}.
    \end{align*}

By using identities of the classical monic Laguerre polynomials, we obtain that the polynomials $Q_{n}(x) = P_{n}(x)\cdot \mathcal{V}$ are given by
$$Q_{n}(x) = \begin{psmallmatrix} \ell_{n}^{(\alpha+2)}(x) & a_{1}\big(\ell_{n+1}^{(\alpha+1)}(x) - x\ell_{n}^{(\alpha+2)}(x)\big) & 0 \\ -a_{1} n\ell_{n-1}^{(\alpha+2)} \;& a_{1}^{2}n \,x\,\ell_{n-1}^{(\alpha+2)}(x) + \ell_{n}^{(\alpha+1)}(x) + a_{2}^{2} (n+\alpha+1)(n+\alpha)n x\, \ell_{n-1}^{(\alpha)}(x)\; & -a_{2} (n+\alpha+1)(n+\alpha)n \ell_{n-1}^{(\alpha)}(x) \\ 0 & a_{2}\big(\ell_{n+1}^{(\alpha+1)}(x) - x\ell_{n}^{(\alpha)}(x)\big) & \ell_{n}^{(\alpha)}(x) \end{psmallmatrix}.$$

The sequence $Q_{n}(x)$ is an eigenfunction of the second-order differential operator
$$E= \partial^{2} xI + \partial \begin{psmallmatrix} \alpha+3 - x & a_{1}x & 0 \\ 0 & \alpha + 2 -x & 0 \\ 0 & 3a_{2}x & \alpha+1-x \end{psmallmatrix} + \begin{psmallmatrix} 0 & a_{1}(\alpha+2) & 0 \\ 0 & 1 & 0 \\ 0 & a_{2} (\alpha+2) & 0\end{psmallmatrix},$$
and we also found that  $Q_{n}(x)$ are orthogonal  polynomials with respect to the  weight
$$\widetilde W(x) = e^{-x}x^{\alpha} \begin{pmatrix} a_{1}^{2}x^{3} + x^{2} & a_{1}x^{2} & a_{1}a_{2}x^{3} \\ a_{1}x^{2} & x & a_{2}x^{2} \\ a_{1}a_{2}x^{3} & a_{2}x^{2} & a_{2}^{2}x^{3}+1 \end{pmatrix}.$$

Thus, 
$$ P_{n}(x)\cdot \mathcal{V}=Q_n(x)=A_n \widetilde P_n(x), $$
where $\widetilde P_n(x)$ is the sequence of monic orthogonal polynomials with respect to $\widetilde W$
and $A_n= \left(\begin{smallmatrix}
1     & -a_1(\alpha+n+2) & 0 \\ 0 &  n(n-1)(n-2)a_{2}^{2} + n(n-1)a_{2}^{2}(\alpha+2)+n(a_{2}^{2}(\alpha+2)(\alpha+1)+a_{1}^{2})+1 &  0\\0& -a_2(\alpha+2+3n) &1
\end{smallmatrix}\right)$.
Therefore $\widetilde W$ is a Darboux transformation of  $W= w_{\alpha+2}\oplus w_{\alpha+1}\oplus w_{\alpha}$.

It is worth noting that the Darboux transformation between $W$ and $\widetilde W$ cannot arise from the factorization of any differential operator $D \in \mathcal{D}(W)$ of order less than six. This fact can be verified by using the expressions of the orthogonal polynomials $Q_{n}$ and $P_{n}$.

\subsection{Example 2} 
We consider an example of matrix-valued orthogonal polynomials arising from matrix-valued spherical functions associated with sphere three-dimensional studied first in \cite{PZ16}.
$$W(x) = (1-x^{2})^{\frac{r}{2}-1} \begin{pmatrix}a(x^{2}-1)+r && -rx \\ -rx && (r-a)(x^{2}-1)+r \end{pmatrix}, $$
with $r >0$, and $0 < a < \frac{r}{2}$.

The algebra $\mathcal{D}(W)$ was comprehensively studied in \cite{Z17}. Subsequently, in \cite{CY18}, the authors identified an orthogonal system for $\mathcal{D}(W)$, proving that the algebra is full. As a consequence, $W$ is a Darboux transformation of  $w(x) = (1 - x^2)^{r/2} I$, the direct sum of classical scalar Gegenbauer weights.
Here, we explicitly prove this by giving an explicit factorization of an operator $D \in \mathcal{D}(w)$. Consider the second-order differential operator 
$$D = \begin{pmatrix} (a-r)\delta_{\frac{r}{2},\frac{r}{2}} +(a-r)^{2}(a+1) & 0 \\ 0 & -a\delta_{\frac{r}{2},\frac{r}{2}} - a^{2}(a-r-1) \end{pmatrix} \in \mathcal{D}(w)$$
where $\delta_{\frac{r}{2},\frac{r}{2}} = \partial^{2} (1-x^{2}) + \partial (-x(r+2))$ is the second-order differential operator of the scalar Jacobi weight $(1-x^{2})^{\frac{r}{2}}$. 
We factorize $D$ as $D = \mathcal{V}\mathcal{N}$, where 
\begin{align*}
        \mathcal{V} & = \partial\begin{pmatrix}-1 && -x \\ -x && -1 \end{pmatrix} + \begin{pmatrix}0 && a-r \\ -a && 0 \end{pmatrix}, \displaybreak[0] \\
   \intertext{ and }     \mathcal{N} & = \partial \begin{pmatrix}r-a && -ax \\ (a-r)x && a\end{pmatrix} + \begin{pmatrix}0 && a(a-r-1) \\ (a+1)(a-r) && 0 \end{pmatrix}.
    \end{align*}
Let $J_{n}^{(\frac{r}{2},\frac{r}{2})}(x)$ be the sequence of monic orthogonal polynomials for the scalar weight $(1-x^{2})^{\frac{r}{2}}$. It follows that $Q_{n}(x) = J_{n}^{(\frac{r}{2},\frac{r}{2})}(x)\cdot \mathcal{V}$ is a sequence of orthogonal polynomials for $W$. Thus, $W$ is a Darboux transformation of the direct sum of Gegenbauer scalar weights $w(x) = (1-x^{2})^{\frac{r}{2}}$.

\begin{remark}

With the explicit factorization of the operator $D = \mathcal{V}\mathcal{N}$, along with techniques similar to those discussed in Section \ref{Sect-sph}, we can determine the algebra $\mathcal{D}(W)$ using an alternative and concise argument compared to that presented in \cite{Z17}.

\end{remark}

\subsection{Example 3} The following Hermite-type weight matrix appears in \cite{DG04} as a solution of the Matrix Bochner problem. 
Let $a, \, b \in \mathbb{R} - \{ 0\}$, and $W$ be the weight matrix 
$$W(x) = e^{-x^{2}} \begin{pmatrix}  1 + a^{2} x^{2} + \frac{a^{2}b^{2}x^{4}}{4} && ax + \frac{ab^{2}x^{3}}{2} &&  \frac{abx^{2}}{2} \\ ax + \frac{ab^{2}x^{3}}{2} && b^{2}x^{2}+1 && bx \\ \frac{abx^{2}}{2} &&  bx && 1\end{pmatrix}.$$
The algebra $\mathcal{D}(W)$ contains the second-order $W$-symmetric differential operator $$E = \partial^{2} I + \partial \left(\begin{smallmatrix} -2x & 2a & 0 \\ 0 & -2x & 2b \\ 0 & 0 & -2x \end{smallmatrix}\right) + \left(\begin{smallmatrix} 4 & 0 & ab \\ 0 & -2 & 0 \\ 0 & 0 & 0 \end{smallmatrix}\right).$$

We will see that $W$ is a Darboux transformation of $\widetilde W= w(x)I$, the direct sum of three copies of the scalar Hermite weights $w(x) = e^{-x^{2}}$. 

After some computations, we see that it is not possible to find a convenient factorization of any lower-order operator, so we need to consider eighth-order operators. From  Theorem \ref{algdirect} we know that any differential operator in the algebra 
$\mathcal{D}(\widetilde W)$ is  a matrix polynomial in $\delta = \partial^{2} + \partial(-2x)$, the classical   (second-order) Hermite  operator, i.e., $$\mathcal{D}(\widetilde W) = \operatorname{Mat}_{3}(\mathbb C)[\delta].$$

We take the eighth-order differential operator $D \in \mathcal{D}(\widetilde W)$ 
\begin{align*}
        D & = \delta^{4} \begin{psmallmatrix} 1 & 0 & 0 \\ 0 & 1 & 0 \\ 0 & 0 & 1 \end{psmallmatrix} + \delta^{3} \begin{psmallmatrix} \frac{2(3a^{2}b^{2}+8a^{2}-8b^{2})}{a^{2}b^{2}} & 0 & -\frac{16}{ab}  \\ 0 & \frac{4(a^{2}b^{2}+4a^{2}+4b^{2})}{a^{2}b^{2}} & 0 \\ -\frac{16}{ab} & 0 & -\frac{2(3a^{2}b^{2}+8a^{2}+32)}{a^{2}b^{2}} \end{psmallmatrix} \displaybreak[0]\\ 
        & \quad + \delta^{2} \begin{psmallmatrix} \frac{4(3a^{2}b^{4}+16a^{2}b^{2}-16b^{4}+16a^{2}-32b^{2})}{a^{2}b^{4}} & 0 & \frac{16(a^{2}b^{2}+4a^{2}+8b^{2}+16)}{a^{3}b^{3}} \\ 0 & \frac{4(a^{2}b^{2}+8a^{2}+16b^{2}+96)}{a^{2}b^{2}} & 0 \\\frac{16(a^{2}b^{2}+4a^{2}+8b^{2}+16)}{a^{3}b^{3}} & 0 & \frac{8(a^{4}b^{4}+4a^{4}b^{2}+64a^{2}b^{2}+128a^{2}+128)}{a^{4}b^{4}}\end{psmallmatrix} \displaybreak[0]\\
        & \quad + \delta \begin{psmallmatrix} \frac{8(a^{2}b^{4}+8a^{2}b^{2}-8b^{4}+16a^{2}-80b^{2}-128)}{a^{2}b^{4}} & 0 & -\frac{32(3a^{2}b^{4}+12a^{2}b^{2}-8b^{4}+48b^{2}+128)}{a^{3}b^{5}} \\ 0 & \frac{-64(a^{2}b^{4}+12a^{2}b^{2}+16a^{2}+16b^{2})}{a^{4}b^{4}} & 0 \\ -\frac{32(3a^{2}b^{4}+12a^{2}b^{2}-8b^{4}+48b^{2}+128)}{a^{3}b^{5}} & 0 & \frac{-128(3a^{2}b^{4}+16a^{2}b^{2}+8b^{4}+80b^{2}+128)}{a^{4}b^{6}}\end{psmallmatrix} \displaybreak[0] \\
        & \quad + \begin{psmallmatrix} \frac{256(3a^{2}b^{4}+16a^{2}b^{2}+16a^{2}+16b^{2})}{a^{4}b^{6}} & 0 &  \frac{1024(a^{2}b^{2}+4a^{2}+16)}{a^{5}b^{5}} \\ 
        0 & \frac{2048(b^{2}+2)}{a^{4}b^{4}} & 0 \\ \frac{1024(a^{2}b^{2}+4a^{2}+16)}{a^{5}b^{5}} & 0 & \frac{4096(3a^{2}b^{2}+8a^{2}+16)}{a^{6}b^{6}}\end{psmallmatrix}.
    \end{align*}

    \smallskip
We factorize the operator $D$ as $D = \mathcal{V}\mathcal{N}$, where

\begin{align*}
        \mathcal{V} & = \partial^{4} \begin{psmallmatrix}  1  && -ax && \frac{abx^{2}}{2} \\  0 && 1 && -bx \\ 0 && 0 &&  1\end{psmallmatrix} + \partial^{3} \begin{psmallmatrix}  -2x &&  2ax^{2}-\frac{4}{a} &&  -abx^{3}+\frac{4bx}{a} \\  -\frac{4}{a} &&  0 && 2bx^{2}-\frac{4}{b} \\ 0 && -\frac{4}{b} && -2x\end{psmallmatrix} \\
        & \quad + \partial^{2} \begin{psmallmatrix} \frac{-2(b^{2}-4)}{b^{2}} & \frac{2(a^{2}b^{2}-4a^{2}+4b^{2})x}{ab^{2}} & - \frac{(a^{2}b^{2} - 4a^{2}+8b^{2})x^{2}}{ab} \\ \frac{8x}{a} & -4x^{2}-6 & \frac{2(3b^{2}+8)x}{b} \\ 0 & \frac{16x}{b} & \frac{-2(2a^{2}b^{2}x^{2}+3a^{2}b^{2}+16)}{a^{2}b^{2}}\end{psmallmatrix} \\ 
        & \quad + \partial \begin{psmallmatrix} 0 &&  0 &&  \frac{16x}{ab} \\ \frac{(16b^{2}+32)}{ab^{2}} && -\frac{(8b^{2}+32)x}{b^{2}} && \frac{(8a^{2}+32)}{a^{2}b} \\ -\frac{16x}{ab} &&  \frac{(8a^{2}b^{2}+32b^{2}+128)}{a^{2}b^{3}} && \frac{(4a^{2}-32)x}{a^{2}} \end{psmallmatrix}  + \begin{psmallmatrix}  -\frac{64}{a^{2}b^{2}} && 0 &&  -\frac{(-32b^{2}-64)}{ab^{3}} \\ 0 && -\frac{64}{a^{2}b^{2}} && 0 \\  -\frac{256}{a^{3}b^{3}} && 0 && 0\end{psmallmatrix}, \displaybreak[0] \\
\intertext {and }        \mathcal{N} & = \partial^{4}  \begin{psmallmatrix} 1 && ax &&  \frac{abx^{2}}{2} \\ 0 &&  1 &&  bx \\ 0 && 0 && 1\end{psmallmatrix} + \partial^{3} \begin{psmallmatrix} -2x && -\frac{(2a^{2}x^{2}-4a^{2}-4)}{a} && -\frac{ax(b^{2}x^{2}-4b^{2}-4)}{b} \\ \frac{4}{a} && 0 &&  -\frac{(2b^{2}x^{2}-4b^{2}-4)}{b} \\ 0 && \frac{4}{b} && -2x \end{psmallmatrix} \\ 
        & \quad + \partial^{2} \begin{psmallmatrix}  -\frac{(4b^{2}x^{2}-4b^{2}-8)}{b^{2}} && -\frac{2x(3a^{2}+8)}{a} && -\frac{(6a^{2}b^{2}x^{2}-6a^{2}b^{2}+8a^{2}x^{2}-12a^{2}+16x^{2})}{ab} \\ -\frac{16x}{a} &&  -4x^{2}+6 && -\frac{2x(3a^{2}b^{2}+4a^{2}+16)}{a^{2}b} \\  0 && -\frac{8x}{b} && -\frac{32}{a^{2}b^{2}}\end{psmallmatrix} \\
        & \quad + \partial \begin{psmallmatrix} -\frac{4x(5b^{2}+8)}{b^{2}} && -\frac{(24b^{2}+32)}{ab^{2}} && -\frac{2x(3a^{2}b^{4}+8a^{2}b^{2}+40b^{2}+64)}{ab^{3}} \\  -\frac{8}{a} && -\frac{8x(a^{2}+4)}{a^{2}} && -\frac{(96b^{2}+128)}{a^{2}b^{3}} \\ -\frac{16x}{ab} &&  -\frac{32}{a^{2}b} && 0 \end{psmallmatrix} \\
        & \quad + \begin{psmallmatrix}  -\frac{(4a^{2}b^{2}+16a^{2}+64)}{a^{2}b^{2}} && 0 &&  -\frac{(48a^{2}b^{2}+128a^{2}+256)}{a^{3}b^{3}} \\ 0 && -\frac{(32b^{2}+64)}{a^{2}b^{2}} && 0 \\ \frac{(16b^{2}+64)}{ab^{3}} && 0 && -\frac{64}{a^{2}b^{2}} \end{psmallmatrix}.
    \end{align*}

    \smallskip
    
It follows that the operators $\mathcal V$ and $\mathcal N$ are degree-preserving operators and the sequence $Q_{n}(x) = h_{n}(x) \cdot \mathcal{V}$ is a sequence of orthogonal polynomials for the weight matrix $W$, where $h_{n}(x)$ is the sequence of monic orthogonal polynomials for the classical scalar Hermite weight. Thus, $W(x)$ is a Darboux transformation of the direct sum $\widetilde W(x) = e^{-x^{2}}I$.

\

\bibliographystyle{amsplain} % Estilo de bibliografía
\bibliography{referencias} % Archivo de bibliografía (sin la extensión

\end{document}